\documentclass[11pt]{article}

\usepackage[margin=2.6cm]{geometry}
\usepackage{amsmath,amssymb,amsthm}
\usepackage{tikz}
\usepackage[colorlinks=true,linkcolor=blue,citecolor=blue,urlcolor=blue]{hyperref}

\numberwithin{equation}{section}

\newtheorem{theorem}{Theorem}[section]
\newtheorem{lemma}[theorem]{Lemma}
\newtheorem{proposition}[theorem]{Proposition}

\theoremstyle{definition}
\newtheorem{definition}[theorem]{Definition}
\newtheorem{example}[theorem]{Example}
\theoremstyle{remark}

\title{Global entropy solutions of the Cauchy problem for the Burgers equation with time-variable coefficient and gradually linear damping}

\author{Richard De la cruz\thanks{School of Mathematics and Statistics,
Universidad Pedag\'ogica y Tecnol\'ogica de Colombia, Avenida Central del Norte 39--115,
150003 Tunja, Boyac\'a, Colombia. E-mail: \texttt{richard.delacruz@uptc.edu.co}
(corresponding author).}
\and
\'Oscar Fernando Rojas-Matamoros\thanks{School of Postgraduate Studies, Faculty of Science,
Universidad Pedag\'ogica y Tecnol\'ogica de Colombia, Avenida Central del Norte 39--115,
150003 Tunja, Boyac\'a, Colombia. E-mail: \texttt{oscar.rojas08@uptc.edu.co}}}

\date{August 26, 2026}

\begin{document}

\maketitle

\begin{abstract}
We study the Cauchy problem for the inviscid Burgers equation with a time-variable coefficient and a gradually linear damping term, $u_t+\alpha(t)(u^2/2)_x=-\sigma(t)u$. We introduce a viscous regularization whose viscosity coefficient depends on time through both $\alpha$ and $\sigma$, and we show that, among viscosities depending only on time, it is the only one, up to a constant factor, for which the Hopf--Cole transformation linearizes the problem without imposing any relation between the two coefficients. This yields explicit viscous solutions and, in the vanishing viscosity limit, an explicit representation of a global weak solution of Lax--Oleinik type, for locally bounded initial data whose primitive grows subquadratically. The limit satisfies a time-dependent one-sided Lipschitz condition. For bounded initial data it is the unique entropy solution, and its amplitude is bounded by $\|u_0\|_{L^\infty}e^{-\int_0^t\sigma(r)\,dr}$. The analysis shows that the competition between nonlinear steepening and damping is governed by a single quantity built from both coefficients, the total displacement of a characteristic per unit initial amplitude. For smooth initial data, shocks form if and only if the maximal compression of the initial datum exceeds the reciprocal of this quantity. As applications, we show that geometric spreading delays but never prevents the breaking of cylindrical and spherical acoustic waves, and we solve the Riemann problem explicitly, showing that shock and rarefaction waves remain confined to a bounded region of space when this quantity is finite.
\end{abstract}

\noindent\textbf{Keywords:} Burgers equation; time-dependent coefficients; time-dependent damping; vanishing viscosity; Hopf--Cole transformation; entropy solutions; Lax--Oleinik formula; shock formation.

\smallskip
\noindent\textbf{MSC 2020:} 35L65, 35L67, 35D30, 35B40.

\section{Introduction}\label{sec:intro}

The Burgers equation
\begin{equation}\label{burgers-mu}
u_t+uu_x=\mu u_{xx},\qquad \mu>0,
\end{equation}
is one of the simplest nonlinear dissipative partial differential equations and arises in a wide variety of physical models. It first appeared in the work of Bateman \cite{Bateman1915} and was studied extensively by Burgers \cite{Burgers1939,Burgers1948} as a model of turbulence. It also describes the formation of large-scale structures in cosmology \cite{Zeldovich1970} and has become a paradigm of nonlinear wave propagation \cite{Whitham1974,BecKhanin2007}. The physical setting that motivates the present work is nonlinear acoustics. Lighthill \cite{Lighthill1956} showed that \eqref{burgers-mu} governs the propagation of plane sound waves of finite amplitude. In suitable variables, $t$ denotes the distance traveled by the wave, $x$ is a retarded time measured in a frame moving with the speed of sound, $u$ is proportional to the acoustic velocity, and $\mu$ is the diffusivity of sound, which accounts for viscosity and heat conduction \cite{Crighton1979}. The nonlinear term describes the steepening of the waveform, which occurs because compressions travel faster than rarefactions, and the diffusive term opposes this steepening. In many situations of interest the acoustic Reynolds number, which measures the strength of the nonlinear effects relative to diffusion, is large, and diffusion is effective only in thin layers in which the waveform changes abruptly. In the limit $\mu\to0^+$ one obtains the inviscid Burgers equation $u_t+uu_x=0$, the simplest genuinely nonlinear scalar conservation law. Its smooth solutions break down in finite time whenever the initial waveform contains a compression. Beyond the breaking time the solution has to be understood as a discontinuous weak solution, whose discontinuities represent shock waves.

The connection between the viscous and the inviscid equations was clarified by Hopf \cite{Hopf1950} and Cole \cite{Cole1951}, who showed independently that the transformation $u=-2\mu(\ln w)_x$ maps \eqref{burgers-mu} onto the heat equation $w_t=\mu w_{xx}$. Hopf used this linearization to compute explicitly the vanishing viscosity limit $\mu\to0^+$, and obtained a global weak solution of the inviscid equation given by a variational formula, nowadays known as the Lax--Oleinik formula, which satisfies a one-sided Lipschitz condition. Since weak solutions of conservation laws are not unique in general, admissibility criteria are needed to single out the physically relevant one. For scalar equations in one space dimension this role is played by the condition~E of Oleinik \cite{Oleinik1957} and by the shock inequalities of Lax \cite{Lax1957}, while Kru\v{z}kov \cite{Kruzkov1970} developed the general theory of entropy solutions of scalar balance laws in several space variables.

Plane waves are, however, an idealization. When sound propagates in a duct whose cross-sectional area $h$ varies slowly along the direction of propagation, or when the wave fronts are cylindrical or spherical, the amplitude of the wave also changes because of geometric spreading, and \eqref{burgers-mu} is replaced by the generalized Burgers equation
\begin{equation}\label{gbe}
u_t+uu_x+\frac{h'(t)}{2h(t)}\,u=\mu u_{xx},
\end{equation}
as explained in \cite{Crighton1979,CrightonScott1979,Sachdev1987}. For cylindrical and spherical waves, $h(t)$ is proportional to $t$ and to $t^2$, respectively, and the linear term takes the form $ju/(2t)$, with $j=1,2$. This term is a linear damping with a time-dependent coefficient, which describes the decrease of the amplitude due to geometric spreading. It can be removed by the change of variable $u=(h(0)/h(t))^{1/2}v$, which transforms \eqref{gbe} into
$$
v_t+\Big(\frac{h(0)}{h(t)}\Big)^{1/2}vv_x=\mu v_{xx},
$$
a Burgers equation with a time-dependent coefficient in front of the nonlinear term. A further change of the time variable then leads to a Burgers equation with a time-dependent viscosity determined by $h$. This last form was proposed by Crighton \cite{Crighton1979} as a model for finite-amplitude sound waves in ducts of variable cross-section, $v$ being an acoustic variable from which the linear effects of the changes of area have been removed. Asymptotic properties and explicit or similarity solutions of these equations were studied in \cite{CrightonScott1979,Scott1981,DoyleEnglefield1990}, and existence, uniqueness and decay rates of classical solutions of the latter equation were obtained in \cite{Zhang2001,WangWarnecke2003,WangZhang2003,ZhangWang2004}. In the limit of large acoustic Reynolds numbers, all these forms lead to inviscid equations of the type $u_t+\alpha(t)uu_x=-\sigma(t)u$, with $\alpha\equiv1$ and $\sigma=h'/(2h)$, or with $\alpha=(h(0)/h)^{1/2}$ and $\sigma\equiv0$. Rao and Yadav \cite{RaoYadav2011} studied the large-time behavior of solutions of the inviscid nonplanar Burgers equation. A natural question, which goes back to the classical theory of nonlinear waves \cite{Whitham1974}, is whether geometric spreading can prevent the formation of shocks.

These considerations lead us to the Cauchy problem
\begin{align}
u_t+\alpha(t)\Big(\frac{u^2}{2}\Big)_x&=-\sigma(t)\,u, && x\in\mathbb{R},\ t>0, \label{main}\\
u(x,0)&=u_0(x), && x\in\mathbb{R}, \label{data}
\end{align}
for the inviscid Burgers equation with a time-variable coefficient $\alpha(t)$ and a gradually linear damping term $\sigma(t)u$, which contains the inviscid limits of both forms of the generalized Burgers equation described above and allows us to treat them in a unified way. Throughout the paper we assume the following hypotheses.
\begin{itemize}
\item[(H1)] $\alpha\in C([0,\infty))$ and $\alpha(t)>0$ for all $t\ge0$.
\item[(H2)] $\sigma\in C([0,\infty))$.
\item[(H3)] $u_0\in L^\infty_{\mathrm{loc}}(\mathbb{R})$, and its primitive $U_0(x):=\int_0^xu_0(r)\,dr$ satisfies $U_0(x)=o(x^2)$ as $|x|\to\infty$.
\end{itemize}
No sign condition is imposed on $\sigma$, although the physically relevant situation, in which $\sigma(t)u$ is a dissipation or friction term, corresponds to $\sigma\ge0$. Assumption (H3) is the growth condition introduced by Hopf for the classical Burgers equation, and it allows unbounded initial data, such as $u_0(x)=|x|^{p}$ with $0\le p<1$.

Equation \eqref{main} also arises in other contexts. In open channel hydraulics, the Saint-Venant equations are often written with a Boussinesq correction coefficient $\beta$, which accounts for nonuniform velocity profiles across the channel section, and with a friction term acting on the averaged velocity. In the balance law system $u_t+\beta(t)uu_x=\Gamma(x,t)-\sigma(t)u$, $\rho_t+\beta(t)(\rho u)_x=0$, studied in \cite{SuarezDelacruz2026}, $u$ is the averaged flow velocity, $\rho$ the cross-sectional area and $\Gamma$ an external forcing. In the absence of forcing, the equation for the velocity decouples and coincides with \eqref{main}, with $\alpha=\beta$. In that work the friction coefficient was tied to the Boussinesq coefficient through $\sigma=(\ln\beta)'$, whereas from the modeling point of view the two coefficients are independent. More generally, damping terms whose coefficients decay in time have attracted considerable attention in the study of the compressible Euler equations and $p$-systems, where the decay rate $(1+t)^{-1}$ turns out to be critical \cite{LiLiMeiZhang2017,GengLinMei2020}. Moreover, the interaction between the nonlinearity of the flux and source terms in scalar equations has been studied, among others, by Murray \cite{Murray1970a,Murray1970b}, who analyzed the decay of discontinuous solutions of perturbed first order wave equations, by Natalini and Tesei \cite{NataliniTesei1992}, who obtained sharp estimates on the support and the decay of entropy solutions of conservation laws with nonlinear absorption, and by Joseph and Sachdev \cite{JosephSachdev1994}. A general account of balance laws with dissipative sources is given in \cite{Dafermos2016}. In contrast with these works, the source term in \eqref{main} is linear, but its coefficient, as well as the coefficient $\alpha$, depends on time in an arbitrary way.

The method of vanishing viscosity requires the choice of a viscous regularization, and for nonautonomous problems this choice is not innocuous. 
In the study of the Riemann problem, Tup\v{c}iev \cite{Tupciev1973} and Dafermos \cite{Dafermos1973} introduced the regularization $u_t+f(u)_x=\varepsilon tu_{xx}$, which preserves the invariance of the Riemann problem under the dilations of the $(x,t)$-plane. For the Burgers equation with Riemann data, Slemrod \cite{Slemrod1995} showed that the difference between the solutions of this regularization and of the classical one $u_t+f(u)_x=\varepsilon u_{xx}$ tends to zero in suitable integral norms as $\varepsilon\to0^+$ for $0<t\le1$. Wang and Zhang \cite{WangZhang2000} introduced the viscosity $\varepsilon e^tu_{xx}$ and extended this comparison to arbitrary bounded time intervals. For systems with damping, self-similar solutions are not available, and regularizations $u_t+f(u)_x=\varepsilon g(t)u_{xx}$ with a nonlinear function $g$ adapted to the damping were introduced in \cite{Delacruz2020}. Regularizations of this type were used in \cite{DelacruzLuWang2024} to solve the Riemann problem for $u_t+k\alpha(t)uu_x=-\lambda\sigma(t)u$ by means of a generalized similarity variable, in \cite{DelacruzNeves2025} for a class of nonautonomous quasilinear systems with time-gradually-degenerate damping, and in \cite{SuarezDelacruz2026} for the Cauchy problem of the Saint-Venant type system mentioned above by means of a generalized Hopf--Cole transformation. On the other hand, viscous Burgers equations with variable coefficients, and Hopf--Cole type transformations for them, have been studied in \cite{Chu1965,Sachdev1978,MalomedShrira1991,Xuetal2007,VaganKumaran2008,BuyukasikPashaev2013,WangZhangLiXin2020}. In particular, Wang, Zhang, Li and Xin \cite{WangZhangLiXin2020} constructed explicit solutions of $u_t+\alpha(t)uu_x+\sigma(t)u=\nu(t)u_{xx}$ under suitable relations among the coefficients. The generalized similarity variable is specific to Riemann data, and the generalized Hopf--Cole transformation requires $\sigma=(\ln\beta)'$. The purpose of the present paper is to solve the Cauchy problem \eqref{main}--\eqref{data} by the vanishing viscosity method for general initial data and independent coefficients. An introductory account of scalar conservation laws with linear damping is given in \cite{DelacruzOchoaRojas2026}.

Our starting point is the choice of the regularization. We set
\begin{equation}\label{SaA}
S(t):=\int_0^t\sigma(r)\,dr,\qquad a(t):=\alpha(t)e^{-S(t)},\qquad A(t):=\int_0^ta(r)\,dr,
\end{equation}
and regularize \eqref{main} by
\begin{equation}\label{viscous}
u_t+\alpha(t)\Big(\frac{u^2}{2}\Big)_x+\sigma(t)u=\frac{\varepsilon}{2}\,a(t)\,u_{xx},\qquad \varepsilon>0 .
\end{equation}
The reason for this choice is the following elementary observation, which is the analogue, for \eqref{viscous}, of the reduction of \eqref{gbe} described above. A function $u$ solves \eqref{viscous} if and only if
\begin{equation}\label{transf}
u(x,t)=e^{-S(t)}\,v\big(x,A(t)\big),
\end{equation}
where $v=v(x,\tau)$ solves the viscous Burgers equation $v_\tau+vv_x=\frac{\varepsilon}{2}v_{xx}$ with constant viscosity. The factor $e^{-S}$ removes the damping, a change of dependent variable already used for damped equations in \cite{Delacruz2020}, and the new time $\tau=A(t)$ absorbs the coefficient $\alpha$. The viscosity in \eqref{viscous} is the one that both operations transform into a constant, so that the classical Hopf--Cole transformation applies. We show that, if the viscosity $\varepsilon u_{xx}$ is used instead, the Hopf--Cole transformation linearizes the problem only when $\sigma=\alpha'/\alpha$, which is exactly the structural restriction $\sigma=(\ln\beta)'$ recalled above. For $k=\lambda=1$, the viscosity used in \cite{DelacruzLuWang2024} is $\varepsilon a(t)A(t)u_{xx}$, which the same transformation maps onto the Dafermos regularization $\varepsilon\tau v_{xx}$. Our choice corresponds instead to the constant viscosity of Hopf, which is the natural one for the Cauchy problem.

The quantities $e^{-S(t)}$ and $A(t)$ have a simple meaning, which is best seen along characteristics. If $u$ is a smooth solution of \eqref{main} and $x(t)$ is a characteristic curve, that is, $x'(t)=\alpha(t)u(x(t),t)$, then $\frac{d}{dt}u(x(t),t)=-\sigma(t)u(x(t),t)$, and therefore
\begin{equation}\label{char-intro}
u\big(x(t),t\big)=e^{-S(t)}\,u_0(\xi),\qquad x(t)=\xi+A(t)\,u_0(\xi),\qquad \xi=x(0).
\end{equation}
Thus $e^{-S(t)}$ is the factor by which the damping has reduced the amplitude up to time $t$, and $A(t)$ is the distance traveled up to time $t$ by a characteristic that carries a unit initial amplitude. Since $A$ is increasing, it has a limit
\begin{equation}\label{Ainf}
A_\infty:=\lim_{t\to\infty}A(t)=\int_0^\infty\alpha(t)\,e^{-\int_0^t\sigma(r)\,dr}\,dt\in(0,\infty],
\end{equation}
which is the total distance that such a characteristic can ever travel. We call $A_\infty$ the \emph{maximal displacement} of the characteristics. Since the displacement $A(t)u_0(\xi)$ is proportional to the initial amplitude, $A(t)$ measures a displacement in the variable $x$ per unit of amplitude rather than a length (in the acoustic interpretation above, where $x$ is a retarded time, $A(t)$ is a retarded time per unit amplitude). When the damping is strong enough, or the coefficient $\alpha$ decays fast enough, $A_\infty$ is finite. In this case the characteristics come to a halt after traveling a finite distance, in the same way as a body moving against friction stops after a finite distance. From the mathematical point of view, $A(t)$ is also the time variable of the classical Burgers equation to which \eqref{main} is reduced by \eqref{transf}, so that the whole evolution of \eqref{main} on $[0,\infty)$ corresponds to the evolution of the Burgers equation on the interval $[0,A_\infty)$. Our results show that $A_\infty$ is the single quantity that decides the long-time behavior of the solutions.

We now describe the main contributions of this paper. The first one is a viscous regularization adapted to independent coefficients. Among all viscosities of the form $\nu(t)u_{xx}$, the viscosity in \eqref{viscous} is, up to a constant factor, the only one for which the change of variables \eqref{transf} combined with the Hopf-Cole transformation yields a linear problem. Previous explicit approaches to this equation either required a relation between the coefficients, such as $\sigma=(\ln\beta)'$ or the conditions of \cite{WangZhangLiXin2020}, or used a Dafermos-type regularization designed to Riemann problem. Our regularization imposes no relation between $\alpha$ and $\sigma$.

The second contribution is the construction of explicit global solutions for general, possibly unbounded, data. For every $\varepsilon>0$, the viscous problem \eqref{viscous}, \eqref{data} has an explicit smooth solution, and as $\varepsilon\to0^+$ these solutions converge almost everywhere to the function
\begin{equation}\label{LO-intro}
u(x,t)=e^{-S(t)}\,\frac{x-y(x,t)}{A(t)},\qquad
y(x,t)\in\operatorname*{arg\,min}_{y\in\mathbb{R}}\Big\{\frac{(x-y)^2}{2A(t)}+U_0(y)\Big\},
\end{equation}
which is a global weak solution of \eqref{main}--\eqref{data}. This is a Lax--Oleinik formula for \eqref{main}. It is Hopf's formula for the Burgers equation, evaluated at the time $A(t)$ and multiplied by the damping factor $e^{-S(t)}$. The solution satisfies the one-sided Lipschitz condition
$$
u(x_2,t)-u(x_1,t)\le e^{-S(t)}\,\frac{x_2-x_1}{A(t)},\qquad x_1<x_2,\ t>0,
$$
so that all its discontinuities are admissible shocks. The result holds for locally bounded data under Hopf's growth condition (H3), and therefore extends Hopf's classical theory to \eqref{main}. For general coefficients $\alpha$ and $\sigma$, the explicit results previously available for \eqref{main} concerned only Riemann data or the case $\sigma=(\ln\beta)'$, and the large-time behavior in the special case $\alpha\equiv1$, $\sigma(t)=j/(2(1+t))$ was studied in \cite{RaoYadav2011}.

The third contribution concerns uniqueness and decay. For bounded initial data, the solution \eqref{LO-intro} is the unique entropy solution, it is bounded by $\|u_0\|_{L^\infty}e^{-S(t)}$, and for integrable data it decays like $e^{-S(t)}A(t)^{-1/2}$. We give a self-contained proof of uniqueness by a duality argument in the spirit of Oleinik, carried out on the possibly bounded interval $[0,A_\infty)$ produced by the change of variables \eqref{transf}. In particular, the Riemann solutions constructed in \cite{DelacruzLuWang2024} are unique in this class.

The fourth contribution is a sharp criterion for shock formation. Let $u_0$ be smooth, and let $-\kappa<0$ be the most negative slope of $u_0$. We prove that the solution remains smooth for all times if and only if $\kappa A_\infty\le1$, and that otherwise its gradient blows up at the explicit time $t^*=A^{-1}(1/\kappa)$. The meaning of this criterion is transparent from \eqref{char-intro}. Two characteristics issued from points where $u_0$ decreases with slope close to $-\kappa$ approach each other, and they meet when $A(t)$ reaches $1/\kappa$. Hence a shock forms if and only if the characteristics can travel far enough. For the generalized Burgers equation \eqref{gbe} with $h(t)$ proportional to $(t_0+t)^j$, $t_0>0$, the damping coefficient is $\sigma(t)=j/(2(t_0+t))$, and the spherical case $j=2$ is critical. Cylindrical and spherical waves always break, the breaking time of spherical waves being exponentially large with respect to $1/\kappa$, whereas for a faster growth of the area, as well as in exponential horns, waveforms with a small compression never break. This is reminiscent of the critical role of damping coefficients decaying like $(1+t)^{-1}$ for the compressible Euler equations \cite{GengLinMei2020}. The same criterion shows that shock formation can be prevented by a coefficient $\alpha$ that decays fast enough, even in the absence of damping, and that it always occurs when $\sigma=\alpha'/\alpha$, since the damping and the coefficient $\alpha$ then compensate each other exactly.

Finally, we solve the Riemann problem explicitly and show that, when $A_\infty<\infty$, shock waves and rarefaction waves remain confined to a bounded region of space for all times.

Beyond the specific equation, we believe that these results are useful in three respects. First, they provide a rigorous framework for the inviscid limit of the generalized Burgers equations of nonlinear acoustics, with a criterion for wave breaking that holds for arbitrary area laws $h$. Second, they make \eqref{main} a nonautonomous balance law for which the entropy solution, its viscous approximations and their dependence on the coefficients are completely explicit, which makes it a convenient benchmark for numerical schemes for balance laws with time-dependent coefficients. Third, the explicit representation of the velocity is the building block for the study of coupled systems, such as the Saint-Venant type system described above, in which the velocity transports a density and delta shock waves, that is, solutions containing Dirac measures supported on shock curves, are expected to arise.

For bounded data, existence and uniqueness of entropy solutions of \eqref{main} also follow from Kru\v{z}kov's general theory, which is formulated for bounded data. However, that theory provides neither explicit solutions nor quantitative information on how the solution depends on the coefficients. Since the reduction to the Burgers equation is explicit, part of our analysis is necessarily close to Hopf's work. We nevertheless give complete proofs, because the transfer of the weak and entropy formulations through \eqref{transf} involves test functions that are only $C^1$ in time, because the vanishing viscosity limit for unbounded data requires a localization argument, and because the dependence of every statement on $\alpha$ and $\sigma$ through $S$ and $A$ is the main information we wish to extract.

The outline of the paper is as follows. In Section~\ref{sec:main} we introduce the solution concepts and the variational problem, and state the main results precisely. In Section~\ref{sec:reduction} we establish the reduction to the classical Burgers equation at the level of weak, entropy and viscous solutions, and explain the choice of the viscosity. In Sections~\ref{sec:viscous} and~\ref{sec:variational} we construct the viscous solutions and study the minimizers of the variational problem, and in Section~\ref{sec:limit} we carry out the vanishing viscosity limit. Section~\ref{sec:uniqueness} is devoted to uniqueness and decay, Section~\ref{sec:threshold} to the formation of shocks from smooth data and its acoustic interpretation, and Section~\ref{sec:riemann} to the Riemann problem.

\section{Solution concepts and main results}\label{sec:main}

In this section we fix the notation, introduce the notions of weak and entropy solution that will be used throughout, describe the variational problem that underlies the explicit representation \eqref{LO-intro}, and state our main results.

The functions $S$, $a$ and $A$ in \eqref{SaA} encode the coefficients of \eqref{main} in the form in which they enter the analysis. As explained in the Introduction, $e^{-S(t)}$ is the factor by which the damping has reduced the amplitude up to time $t$, and $A(t)$ is the distance traveled by a characteristic carrying a unit initial amplitude. At the same time, $A(t)$ is the time variable of the reduced Burgers equation. By (H1) and (H2), $a\in C([0,\infty))$ and $a>0$. Hence $A\in C^1([0,\infty))$, $A(0)=0$ and $A'=a>0$, so that $A:[0,\infty)\to[0,A_\infty)$ is a $C^1$ diffeomorphism, with $A_\infty\in(0,\infty]$ given by \eqref{Ainf}. We denote its inverse by $\theta:[0,A_\infty)\to[0,\infty)$, which satisfies
$$
\theta\in C^1([0,A_\infty)),\qquad \theta'(\tau)=\frac{1}{a(\theta(\tau))}.
$$
For $T>0$ we shall frequently use the constants
\begin{equation}\label{constants}
\bar S_T:=\max_{[0,T]}|S|,\qquad \underline a_T:=\min_{[0,T]}a>0,\qquad \bar a_T:=\max_{[0,T]}a,
\end{equation}
which satisfy $\underline a_Tt\le A(t)\le\bar a_Tt$ for $0\le t\le T$.

The classical Burgers equation in the time variable $\tau$ will be written as
\begin{equation}\label{burgers}
\begin{cases}
v_\tau+\Big(\frac{v^2}{2}\Big)_x=0,\\
v(x,0)=u_0(x),
\end{cases}
\end{equation}
and its viscous regularization with constant viscosity as
\begin{equation}\label{burgers-visc}
v_\tau+\Big(\frac{v^2}{2}\Big)_x=\frac{\varepsilon}{2}\,v_{xx}.
\end{equation}

Solutions of \eqref{main} develop discontinuities in finite time for generic data, so \eqref{main} has to be understood in the sense of distributions. Since the reduction \eqref{transf} maps \eqref{main} on $[0,\infty)$ onto the Burgers equation on the possibly bounded interval $[0,A_\infty)$, we formulate the definitions on a general time interval $(0,T^*)$. Here $L^\infty_{\mathrm{loc}}(\mathbb{R}\times[0,T^*))$ denotes the set of measurable functions that are essentially bounded on $[-R,R]\times(0,T)$ for all $R>0$ and $0<T<T^*$.

\begin{definition}[Weak solution]\label{def:weak}
Let $T^*\in(0,\infty]$. A function $u\in L^\infty_{\mathrm{loc}}(\mathbb{R}\times[0,T^*))$ is a \emph{weak solution of \eqref{main}--\eqref{data} in $\mathbb{R}\times(0,T^*)$} if
\begin{equation}\label{weak}
\int_0^{T^*}\!\!\int_{\mathbb{R}}\Big[u\varphi_t+\alpha(t)\frac{u^2}{2}\varphi_x-\sigma(t)u\varphi\Big]dx\,dt
+\int_{\mathbb{R}} u_0(x)\varphi(x,0)\,dx=0
\end{equation}
for every $\varphi\in C_c^\infty(\mathbb{R}\times[0,T^*))$. When $T^*=\infty$ we simply speak of a weak solution. With $\alpha\equiv1$ and $\sigma\equiv0$, this defines weak solutions of \eqref{burgers} in $\mathbb{R}\times(0,T^*)$.
\end{definition}

Weak solutions are not unique, and an admissibility condition is needed. For the classical Burgers equation, Oleinik's condition~E requires that $v(x+z,\tau)-v(x,\tau)\le Cz/\tau$ for $z>0$. Thus the solution may jump down, but not up, and rarefactions spread at most linearly. For \eqref{main}, the transformation \eqref{transf} multiplies the solution by $e^{\pm S}$ and replaces $t$ by $A(t)$, and both operations distort one-sided bounds by factors that depend on the time interval under consideration. This motivates the following definition, in which the constant is allowed to depend on $T$.

\begin{definition}[Entropy solution]\label{def:entropy}
Assume, in addition to (H1)--(H2), that $u_0\in L^\infty(\mathbb{R})$. A weak solution $u$ of \eqref{main}--\eqref{data} is an \emph{entropy solution} if, for every $T>0$, $u\in L^\infty(\mathbb{R}\times(0,T))$ and there is a constant $C_T\ge0$ such that
\begin{equation}\label{oleinik}
u(x+z,t)-u(x,t)\le C_T\Big(1+\frac1t\Big)z
\quad\text{for a.e. }(x,z,t)\in\mathbb{R}\times(0,\infty)\times(0,T).
\end{equation}
For \eqref{burgers} in $\mathbb{R}\times(0,T^*)$ with $T^*<\infty$, we require \eqref{oleinik} on $\mathbb{R}\times(0,\infty)\times(0,T^*)$ with a single constant.
\end{definition}

For $\alpha\equiv1$ and $\sigma\equiv0$, Definition~\ref{def:entropy} is the notion of entropy solution used by Oleinik \cite{Oleinik1957}, which is also presented in \cite[Section~3.4.3]{Evans2010}. The additive term $C_T$ in $C_T(1+1/t)$ is harmless for the classical equation, but it is convenient here because the factors $e^{\pm S}$ are not identically one.

Hopf's explicit solution of the Burgers equation is obtained by minimizing, for each $(x,\tau)$, the sum of a quadratic transport cost and the primitive of the initial datum. The same variational problem, evaluated at $\tau=A(t)$, describes the solutions of \eqref{main}. For $x,y\in\mathbb{R}$ and $\tau>0$ we set
\begin{equation}\label{G}
G(x,y,\tau):=\frac{(x-y)^2}{2\tau}+U_0(y)
\quad \mbox{ and } \quad
\mathcal M(x,\tau):=\Big\{z\in\mathbb{R}:\ G(x,z,\tau)=\min_{y\in\mathbb{R}}G(x,y,\tau)\Big\}.
\end{equation}
Under (H3), $\mathcal M(x,\tau)$ is a nonempty compact set. Indeed, $y\mapsto G(x,y,\tau)$ is continuous because $U_0$ is locally Lipschitz. Moreover, $(x-y)^2\ge\frac12y^2-x^2$ and $|U_0(y)|\le y^2/(8\tau)$ for $|y|$ large, so that $G(x,y,\tau)\ge y^2/(8\tau)-x^2/(2\tau)\to+\infty$ as $|y|\to\infty$. Hence the minimum is attained and the set of minimizers is closed and bounded. We define the smallest and the largest minimizer,
\begin{equation*}\label{ypm}
y_-(x,\tau):=\min\mathcal M(x,\tau)
\qquad \mbox{ and } \qquad 
y_+(x,\tau):=\max\mathcal M(x,\tau).
\end{equation*}
For each fixed $\tau$, the two coincide except at countably many points $x$, because the minimizers are monotone with respect to $x$. The choice of $y_+$ in the statements below therefore only fixes a right-continuous representative.

We now state our main results. The first one describes the viscous approximations and their limit, and holds for locally bounded initial data.

\begin{theorem}[Existence and vanishing viscosity]\label{thm:existence}
Assume (H1)--(H3).
\begin{itemize}
\item[(i)] For every $\varepsilon>0$ the function
\begin{equation}\label{ueps}
u^\varepsilon(x,t):=e^{-S(t)}\,
\frac{\displaystyle\int_{\mathbb{R}}\frac{x-y}{A(t)}\exp\Big\{-\frac1\varepsilon G\big(x,y,A(t)\big)\Big\}dy}
{\displaystyle\int_{\mathbb{R}}\exp\Big\{-\frac1\varepsilon G\big(x,y,A(t)\big)\Big\}dy},
\qquad x\in\mathbb{R},\ t>0,
\end{equation}
is well defined. The functions $u^\varepsilon$, $u^\varepsilon_x$, $u^\varepsilon_{xx}$ and $u^\varepsilon_t$ are continuous in $\mathbb{R}\times(0,\infty)$, $u^\varepsilon$ satisfies \eqref{viscous} pointwise in $\mathbb{R}\times(0,\infty)$, and $\lim\limits_{t\to0^+}\int_0^xu^\varepsilon(z,t)\,dz=U_0(x)$ for every $x\in\mathbb{R}$.
\item[(ii)] Let $y(x,t):=y_+(x,A(t))$ and
\begin{equation}\label{u}
u(x,t):=e^{-S(t)}\,\frac{x-y(x,t)}{A(t)},\qquad x\in\mathbb{R},\ t>0.
\end{equation}
Then $u$ is a Borel function, and there is a Borel set $\mathcal N\subset\mathbb{R}\times(0,\infty)$ of Lebesgue measure zero, each of whose sections $\mathcal N_t=\{x:(x,t)\in\mathcal N\}$ is countable, such that for every $(x,t)\notin\mathcal N$ the function $y\mapsto G(x,y,A(t))$ has a unique minimizer and $u^\varepsilon(x,t)\to u(x,t)$ as $\varepsilon\to0^+$.
\item[(iii)] For all $R,T>0$ there is $L=L(R,A(T),u_0)$ such that $|u(x,t)|\le Le^{-S(t)}$ for $|x|\le R$ and $0<t\le T$. In particular $u\in L^\infty_{\mathrm{loc}}(\mathbb{R}\times[0,\infty))$.
\item[(iv)] $u$ is a weak solution of \eqref{main}--\eqref{data}.
\item[(v)] For every $t>0$ and $x_1<x_2$,
\begin{equation}\label{oleinik-u}
u(x_2,t)-u(x_1,t)\le e^{-S(t)}\,\frac{x_2-x_1}{A(t)} .
\end{equation}
Consequently, for every $t>0$, $u(\cdot,t)$ is of locally bounded variation and right-continuous, its discontinuity points form a countable set, and $u(x-,t)>u(x+,t)$ at each of them.
\end{itemize}
\end{theorem}

Property (v) states that every discontinuity of $u$ is a downward jump, as for admissible shocks of the classical Burgers equation. It is the property that singles out the solution among all weak solutions when the data are bounded.

\begin{theorem}[Uniqueness and decay for bounded data]\label{thm:uniqueness}
Assume (H1)--(H2) and $u_0\in L^\infty(\mathbb{R})$. Then (H3) holds, and the function $u$ in \eqref{u} has the following properties.
\begin{itemize}
\item[(i)] $|u(x,t)|\le\|u_0\|_{L^\infty}\,e^{-S(t)}$ for all $x\in\mathbb{R}$ and $t>0$.
\item[(ii)] $u$ is an entropy solution of \eqref{main}--\eqref{data}, and every entropy solution $\hat u$ satisfies $\hat u=u$ a.e.\ in $\mathbb{R}\times(0,\infty)$.
\item[(iii)] If moreover $u_0\in L^1(\mathbb{R})$, then
\begin{equation}\label{decay-L1}
|u(x,t)|\le e^{-S(t)}\Big(\frac{2\|u_0\|_{L^1}}{A(t)}\Big)^{1/2}
\qquad\text{for all }x\in\mathbb{R},\ t>0 .
\end{equation}
\end{itemize}
\end{theorem}

The two decay mechanisms in Theorem~\ref{thm:uniqueness} have different origins. If $S(t)\to+\infty$, part (i) gives uniform decay at the rate $e^{-S(t)}$, independently of $\alpha$. If $\sigma\equiv0$ and $A_\infty=\int_0^\infty\alpha=\infty$, then \eqref{decay-L1} gives the decay $A(t)^{-1/2}$, which reduces to the classical rate $t^{-1/2}$ when $\alpha\equiv1$. If $A_\infty<\infty$, the factor $A(t)^{-1/2}$ stays bounded away from zero, and the only decay mechanism is $e^{-S(t)}$.

The next result shows that the maximal displacement $A_\infty$ determines whether smooth solutions remain smooth.

\begin{theorem}[Global smoothness versus shock formation]\label{thm:threshold}
Assume (H1)--(H2), let $u_0\in C^1(\mathbb{R})$ with $u_0,u_0'$ bounded, and set $\kappa:=\max\{0,-\inf\limits_{\mathbb{R}} u_0'\}$.
\begin{itemize}
\item[(i)] If $\kappa A_\infty\le1$ (which, when $A_\infty=\infty$, means $\kappa=0$), then $u\in C^1(\mathbb{R}\times[0,\infty))$ is a classical solution of \eqref{main}--\eqref{data}, given by
\begin{equation}\label{characteristics}
u(x,t)=e^{-S(t)}u_0(\xi),\qquad x=\xi+A(t)u_0(\xi),
\end{equation}
and, in $\mathbb{R}\times[0,\infty)$,
$$
-\frac{\kappa\,e^{-S(t)}}{1-\kappa A(t)}\le u_x(x,t)\le e^{-S(t)}\|u_0'\|_{L^\infty}.
$$
In particular, $u(\cdot,t)$ is Lipschitz continuous for every $t\ge0$.
\item[(ii)] If $\kappa A_\infty>1$, then $t^*:=\theta(1/\kappa)\in(0,\infty)$, $u\in C^1(\mathbb{R}\times[0,t^*))$ is given by \eqref{characteristics} for $0\le t<t^*$, and
\begin{equation}\label{blowup}
\sup_{x\in\mathbb{R}}\big(-u_x(x,t)\big)=\frac{\kappa\,e^{-S(t)}}{1-\kappa A(t)}\longrightarrow+\infty
\qquad\text{as }t\to(t^*)^- .
\end{equation}
\end{itemize}
\end{theorem}

Finally, for piecewise constant data with a single jump, the solution can be computed in closed form.

\begin{theorem}[Riemann problem]\label{thm:riemann}
Assume (H1)--(H2) and let $u_0(x)=u_l$ for $x<0$ and $u_0(x)=u_r$ for $x>0$, with constants $u_l\ne u_r$. Then the unique entropy solution of \eqref{main}--\eqref{data} is
\begin{equation}\label{rarefaction}
u(x,t)=
\begin{cases}
u_l\,e^{-S(t)}, & x<u_lA(t),\\[2pt]
\dfrac{x}{A(t)}\,e^{-S(t)}, & u_lA(t)\le x\le u_rA(t),\\[6pt]
u_r\,e^{-S(t)}, & x>u_rA(t),
\end{cases}
\qquad\text{if }u_l<u_r,
\end{equation}
and
\begin{equation}\label{shock}
u(x,t)=
\begin{cases}
u_l\,e^{-S(t)}, & x<s(t),\\
u_r\,e^{-S(t)}, & x>s(t),
\end{cases}
\qquad s(t):=\frac{u_l+u_r}{2}\,A(t),
\qquad\text{if }u_l>u_r .
\end{equation}
\end{theorem}

\section{Reduction to the classical Burgers equation}\label{sec:reduction}

The transformation \eqref{transf} is the backbone of the paper. For smooth solutions, the equivalence between \eqref{main} and the Burgers equation is a direct application of the chain rule. For weak solutions some care is needed, because test functions for \eqref{main} and for \eqref{burgers} are related by the factor $e^{S(t)}$ and the time change $\tau=A(t)$, and since $S$ and $A$ are only of class $C^1$, the transformed test functions are not smooth. In this section we show that this loss of regularity is harmless, and we transfer weak solutions, one-sided bounds and viscous solutions from one problem to the other. We also explain why the viscosity in \eqref{viscous} is the only one for which the Hopf--Cole transformation is available when $\alpha$ and $\sigma$ are independent.

We first enlarge the class of admissible test functions, so that test functions of class $C^1$ or merely Lipschitz may be used. This will be needed both for the change of variables and for the duality argument that yields uniqueness.

\begin{lemma}[Lipschitz test functions]\label{lem:lipschitz}
Let $T^*\in(0,\infty]$ and let $u\in L^\infty_{\mathrm{loc}}(\mathbb{R}\times[0,T^*))$ satisfy \eqref{weak} for every $\varphi\in C_c^\infty(\mathbb{R}\times[0,T^*))$. Then \eqref{weak} holds for every Lipschitz continuous function $\varphi$ with compact support in $\mathbb{R}\times[0,T^*)$, where $\varphi_t,\varphi_x$ denote the a.e.\ defined derivatives.
\end{lemma}

\begin{proof}
Let $\varphi$ be Lipschitz with $\operatorname{supp}\varphi\subset[-R,R]\times[0,T]$, $T<T^*$, and let $\ell$ be its Lipschitz constant. The even extension $\bar\varphi(x,t):=\varphi(x,|t|)$ is Lipschitz on $\mathbb{R}^2$ with constant $\ell$ and support in $[-R,R]\times[-T,T]$. Let $J_\delta$ be a standard mollifier on $\mathbb{R}^2$ and $\varphi_\delta:=J_\delta*\bar\varphi\in C_c^\infty(\mathbb{R}^2)$. Fix $\delta_0\in(0,1]$ with $T+\delta_0<T^*$. For $0<\delta<\delta_0$, the restriction of $\varphi_\delta$ to $\mathbb{R}\times[0,T^*)$ belongs to $C_c^\infty(\mathbb{R}\times[0,T^*))$ and is supported in the compact set $[-R-1,R+1]\times[0,T+\delta_0]$. Moreover $\|\varphi_\delta-\bar\varphi\|_{L^\infty}\le\ell\delta$, $|\nabla\varphi_\delta|\le\ell$, and $\nabla\varphi_\delta=J_\delta*\nabla\bar\varphi\to\nabla\bar\varphi$ at every Lebesgue point of $\nabla\bar\varphi$, hence a.e.\ in $\mathbb{R}^2$. Since $u$, $u^2$, $\alpha$ and $\sigma$ are bounded on this set and $u_0\in L^1(-R-1,R+1)$, the dominated convergence theorem allows us to pass to the limit $\delta\to0^+$ in each term of \eqref{weak} written for $\varphi_\delta$.
\end{proof}

We can now transfer weak solutions. Throughout, $u$ and $v$ are related by
\begin{equation}\label{v-from-u}
v(x,\tau):=e^{S(\theta(\tau))}\,u\big(x,\theta(\tau)\big),\qquad (x,\tau)\in\mathbb{R}\times(0,A_\infty),
\end{equation}
which is the inverse of \eqref{transf}.

\begin{lemma}[Reduction of weak solutions]\label{lem:weak-reduction}
Let $u:\mathbb{R}\times(0,\infty)\to\mathbb{R}$ be measurable and let $v$ be given by \eqref{v-from-u}. Then $u\in L^\infty_{\mathrm{loc}}(\mathbb{R}\times[0,\infty))$ if and only if $v\in L^\infty_{\mathrm{loc}}(\mathbb{R}\times[0,A_\infty))$, and $u$ is a weak solution of \eqref{main}--\eqref{data} if and only if $v$ is a weak solution of \eqref{burgers} in $\mathbb{R}\times(0,A_\infty)$.
\end{lemma}

\begin{proof}
The map $\Lambda:\mathbb{R}\times[0,\infty)\to\mathbb{R}\times[0,A_\infty)$ defined by $\Lambda(x,t):=(x,A(t))$ is a $C^1$ diffeomorphism. Both $\Lambda$ and $\Lambda^{-1}$ are locally Lipschitz, so they map compact sets onto compact sets and null sets onto null sets. Since $e^{\pm S}$ is bounded on compact time intervals, the first assertion follows.

Let $\psi\in C^1_c(\mathbb{R}\times[0,A_\infty))$ and define
\begin{equation}\label{test-transf}
\varphi(x,t):=e^{S(t)}\,\psi\big(x,A(t)\big).
\end{equation}
If $\operatorname{supp}\psi\subset[-R,R]\times[0,\tau_1]$ with $\tau_1<A_\infty$, then $\varphi\in C^1_c(\mathbb{R}\times[0,\infty))$ with support in $[-R,R]\times[0,\theta(\tau_1)]$. Conversely, every $\varphi\in C^1_c(\mathbb{R}\times[0,\infty))$ is of the form \eqref{test-transf} with
$$
\psi(x,\tau):=e^{-S(\theta(\tau))}\varphi\big(x,\theta(\tau)\big)\in C^1_c(\mathbb{R}\times[0,A_\infty)),
$$
because $\theta\in C^1([0,A_\infty))$ and $S\in C^1([0,\infty))$. By the chain rule,
$$
\varphi_t=\sigma(t)\varphi+e^{S(t)}a(t)\,\psi_\tau\big(x,A(t)\big),\qquad
\varphi_x=e^{S(t)}\psi_x\big(x,A(t)\big).
$$
Since $u(x,t)=e^{-S(t)}v(x,A(t))$, a direct computation gives, for $t>0$,
\begin{align*}
u\varphi_t&=\sigma(t)\,v\psi+a(t)\,v\psi_\tau,\qquad
-\sigma(t)u\varphi=-\sigma(t)\,v\psi,\\
\alpha(t)\frac{u^2}{2}\varphi_x&=\alpha(t)e^{-2S(t)}\frac{v^2}{2}\,e^{S(t)}\psi_x
=a(t)\,\frac{v^2}{2}\psi_x,
\end{align*}
where $v,\psi$ and their derivatives are evaluated at $(x,A(t))$. Hence
\begin{equation}\label{integrand-transf}
u\varphi_t+\alpha(t)\frac{u^2}{2}\varphi_x-\sigma(t)u\varphi
=a(t)\Big[v\psi_\tau+\frac{v^2}{2}\psi_x\Big]\big(x,A(t)\big).
\end{equation}
Both sides are integrable on $\mathbb{R}\times(0,\infty)$ whenever $u$ (equivalently $v$) is locally bounded. Integrating \eqref{integrand-transf} and using the substitution $\tau=A(t)$, $d\tau=a(t)\,dt$, together with $\varphi(x,0)=\psi(x,0)$, we obtain
\begin{equation}\label{functionals}
\int_0^\infty\!\!\int_{\mathbb{R}}\Big[u\varphi_t+\alpha\frac{u^2}{2}\varphi_x-\sigma u\varphi\Big]dx\,dt
+\int_{\mathbb{R}} u_0\varphi(x,0)\,dx
=\int_0^{A_\infty}\!\!\int_{\mathbb{R}}\Big[v\psi_\tau+\frac{v^2}{2}\psi_x\Big]dx\,d\tau
+\int_{\mathbb{R}} u_0\psi(x,0)\,dx .
\end{equation}
If $u$ is a weak solution, Lemma~\ref{lem:lipschitz} shows that the left-hand side of \eqref{functionals} vanishes for every $\varphi\in C^1_c(\mathbb{R}\times[0,\infty))$. Hence the right-hand side vanishes for every $\psi\in C^1_c(\mathbb{R}\times[0,A_\infty))$, in particular for every $\psi\in C_c^\infty(\mathbb{R}\times[0,A_\infty))$. The converse implication is identical.
\end{proof}

The one-sided bound \eqref{oleinik} is transferred in the same way, at the price of constants that depend on the time interval.

\begin{lemma}\label{lem:onesided}
Let $u$ and $v$ be related by \eqref{v-from-u}, let $T>0$, and let $\bar S_T$, $\underline a_T$, $\bar a_T$ be given by \eqref{constants}.
\begin{itemize}
\item[(i)] If $v(x+z,\tau)-v(x,\tau)\le C(1+1/\tau)z$ for a.e.\ $(x,z,\tau)\in\mathbb{R}\times(0,\infty)\times(0,A(T))$, then \eqref{oleinik} holds on $\mathbb{R}\times(0,\infty)\times(0,T)$ with $C_T=Ce^{\bar S_T}\max\{1,1/\underline a_T\}$.
\item[(ii)] If $u$ satisfies \eqref{oleinik} on $\mathbb{R}\times(0,\infty)\times(0,T)$, then $v(x+z,\tau)-v(x,\tau)\le C'(1+1/\tau)z$ for a.e.\ $(x,z,\tau)\in\mathbb{R}\times(0,\infty)\times(0,A(T))$, with $C'=C_Te^{\bar S_T}\max\{1,\bar a_T\}$.
\end{itemize}
\end{lemma}

\begin{proof}
The map $(x,z,t)\mapsto(x,z,A(t))$ preserves null sets. In case (i), for a.e.\ $(x,z,t)$,
$$
u(x+z,t)-u(x,t)=e^{-S(t)}\big[v(x+z,A(t))-v(x,A(t))\big]
\le e^{\bar S_T}C\Big(1+\frac{1}{\underline a_Tt}\Big)z ,
$$
and (i) follows. Case (ii) is analogous, using $1/\theta(\tau)\le\bar a_T/\tau$ for $0<\tau<A(T)$.
\end{proof}

At the viscous level the reduction is a pointwise identity, and it explains the choice of the time-dependent viscosity in \eqref{viscous}.

\begin{lemma}[Reduction of the viscous problem]\label{lem:viscous-reduction}
Let $v^\varepsilon\in C(\mathbb{R}\times(0,\infty))$ have continuous derivatives $v^\varepsilon_x$, $v^\varepsilon_{xx}$, $v^\varepsilon_\tau$, and set $u^\varepsilon(x,t):=e^{-S(t)}v^\varepsilon(x,A(t))$ for $x\in\mathbb{R}$, $t>0$. Then $u^\varepsilon$, $u^\varepsilon_x$, $u^\varepsilon_{xx}$, $u^\varepsilon_t$ are continuous in $\mathbb{R}\times(0,\infty)$, and
$$
u^\varepsilon_t+\alpha u^\varepsilon u^\varepsilon_x+\sigma u^\varepsilon-\frac{\varepsilon}{2}a\,u^\varepsilon_{xx}
=e^{-S(t)}a(t)\Big[v^\varepsilon_\tau+v^\varepsilon v^\varepsilon_x-\frac{\varepsilon}{2}v^\varepsilon_{xx}\Big]\big(x,A(t)\big).
$$
In particular, $u^\varepsilon$ solves \eqref{viscous} in $\mathbb{R}\times(0,\infty)$ if and only if $v^\varepsilon$ solves \eqref{burgers-visc} in $\mathbb{R}\times(0,A_\infty)$.
\end{lemma}

\begin{proof}
The regularity follows from $A,S\in C^1([0,\infty))$. With $v^\varepsilon$ and its derivatives evaluated at $(x,A(t))$, we have $u^\varepsilon_t=-\sigma u^\varepsilon+e^{-S}a\,v^\varepsilon_\tau$, $\alpha u^\varepsilon u^\varepsilon_x=\alpha e^{-2S}v^\varepsilon v^\varepsilon_x=e^{-S}a\,v^\varepsilon v^\varepsilon_x$ and $\frac\varepsilon2 a\,u^\varepsilon_{xx}=e^{-S}a\,\frac\varepsilon2 v^\varepsilon_{xx}$. Since $e^{-S}a>0$, the equivalence follows.
\end{proof}

The factor $a(t)=\alpha(t)e^{-S(t)}$ in \eqref{viscous} is precisely the one for which the regularized problem is mapped onto a viscous Burgers equation with constant viscosity. To see this, assume that $\alpha,\sigma\in C^1([0,\infty))$ and replace the viscous term in \eqref{viscous} by $\nu(t)u_{xx}$, with $\nu\in C^1([0,\infty))$ positive. The computation in the proof of Lemma~\ref{lem:viscous-reduction} shows that the same transformation leads to $v_\tau+vv_x=\bar\nu(\tau)v_{xx}$ with $\bar\nu(\tau)=\nu(\theta(\tau))/a(\theta(\tau))$. Writing $v=V_x$ and $V=-2c(\tau)\ln w$, the integrated equation $V_\tau+\frac12V_x^2=\bar\nu V_{xx}$ becomes
$$
-2c'\ln w-2c\,\frac{w_\tau}{w}+2c(c-\bar\nu)\frac{w_x^2}{w^2}=-2c\bar\nu\frac{w_{xx}}{w}.
$$
The quadratic term disappears only if $c=\bar\nu$, and then $w_\tau=\bar\nu w_{xx}-(\bar\nu'/\bar\nu)w\ln w$. Thus the Hopf--Cole transformation linearizes the problem only when $\bar\nu$ is constant, that is, when $\nu$ is a constant multiple of $a$. For the constant viscosity $\nu\equiv\varepsilon$ this means $a\equiv a(0)$, that is, $\sigma=\alpha'/\alpha$, which is exactly the structural restriction $\sigma=(\ln\beta)'$ imposed in \cite{SuarezDelacruz2026}.

Two particular cases illustrate the choice of the viscosity. When $\alpha\in C^1([0,\infty))$ and $\sigma=\alpha'/\alpha$, we have $e^{S}=\alpha/\alpha(0)$, $a\equiv\alpha(0)$ and $A(t)=\alpha(0)t$, so that \eqref{viscous} has constant viscosity. With $\alpha(0)=1$ it is the regularization $u_t+\beta uu_x=\Gamma-\sigma u+\frac\varepsilon2u_{xx}$ of \cite{SuarezDelacruz2026}, where $\beta(0)=1$ and $\sigma=(\ln\beta)'$, with vanishing source term $\Gamma$. When $\sigma\equiv0$, \eqref{viscous} reads $u_t+\alpha(t)uu_x=\frac\varepsilon2\alpha(t)u_{xx}$.

\section{Viscous solutions by the Hopf--Cole transformation}\label{sec:viscous}

By Lemma~\ref{lem:viscous-reduction}, constructing solutions of \eqref{viscous} amounts to constructing solutions of the viscous Burgers equation \eqref{burgers-visc} with constant viscosity. Following Hopf, we write $v^\varepsilon=-\varepsilon(\ln w)_x$, where $w$ solves the heat equation with the initial datum $e^{-U_0/\varepsilon}$. Since $U_0$ may grow almost quadratically, this initial datum may grow like $e^{\eta x^2}$ for every $\eta>0$, and some care is needed to show that the heat kernel representation converges, is smooth, and attains the initial datum. In this section and in the next one only assumption (H3) is used, and the time variable is $\tau\in(0,\infty)$. The coefficients $\alpha$ and $\sigma$ enter only through the change of variables \eqref{transf}.

The following elementary consequence of (H3) will be used repeatedly.

\begin{lemma}\label{lem:U0-growth}
Under (H3), for every $\eta>0$ there is $K_\eta\ge0$ such that $|U_0(y)|\le\eta y^2+K_\eta$ for all $y\in\mathbb{R}$.
\end{lemma}

\begin{proof}
By (H3) there is $R_\eta>0$ with $|U_0(y)|\le\eta y^2$ for $|y|\ge R_\eta$. We take $K_\eta:=\max\limits_{|y|\le R_\eta}|U_0(y)|$, which is finite because $U_0$ is continuous.
\end{proof}

\begin{lemma}\label{lem:heat}
Assume (H3), let $\varepsilon>0$ and $w_0:=e^{-U_0/\varepsilon}$. For $x\in\mathbb{R}$ and $\tau>0$ let
\begin{equation}\label{heat-sol}
w(x,\tau):=\int_{\mathbb{R}}\Phi_\varepsilon(x-y,\tau)\,w_0(y)\,dy
\qquad \mbox{ and } \qquad
\Phi_\varepsilon(z,\tau):=\frac{1}{\sqrt{2\pi\varepsilon\tau}}\exp\Big\{-\frac{z^2}{2\varepsilon\tau}\Big\}.
\end{equation}
Then the integral converges absolutely, $w\in C^\infty(\mathbb{R}\times(0,\infty))$, $w>0$, $w_\tau=\frac\varepsilon2w_{xx}$ in $\mathbb{R}\times(0,\infty)$, and $w(x,\tau)\to w_0(x_0)$ as $(x,\tau)\to(x_0,0^+)$, for every $x_0\in\mathbb{R}$.
\end{lemma}

\begin{proof}
Let $Q=[-R,R]\times[\tau_1,\tau_2]\subset\mathbb{R}\times(0,\infty)$. Using $(x-y)^2\ge\frac12y^2-x^2$ and Lemma~\ref{lem:U0-growth} with $\eta=1/(8\tau_2)$ we get, for $(x,\tau)\in Q$ and $y\in\mathbb{R}$,
$$
\exp\Big\{-\frac{(x-y)^2}{2\varepsilon\tau}\Big\}w_0(y)
\le\exp\Big\{\frac{R^2}{2\varepsilon\tau_1}+\frac{K_\eta}{\varepsilon}\Big\}
\exp\Big\{-\frac{y^2}{4\varepsilon\tau_2}+\frac{y^2}{8\varepsilon\tau_2}\Big\}
=C_Q\exp\Big\{-\frac{y^2}{8\varepsilon\tau_2}\Big\}.
$$
Every partial derivative of $\Phi_\varepsilon(x-y,\tau)$ with respect to $x$ and $\tau$ is the product of $\Phi_\varepsilon(x-y,\tau)$ and a polynomial in $x-y$ whose coefficients are continuous functions of $\tau>0$. Hence on $Q$ it is bounded by the product of $\Phi_\varepsilon(x-y,\tau)$ and a polynomial in $|y|$ with constant coefficients. Therefore all derivatives of the integrand in \eqref{heat-sol} are dominated on $Q$ by an integrable function of $y$, and differentiation under the integral sign is legitimate. Thus $w\in C^\infty(\mathbb{R}\times(0,\infty))$, and $w$ solves the heat equation because $\Phi_\varepsilon$ does. Finally, $w>0$ because the integrand in \eqref{heat-sol} is positive.

Fix $x_0\in\mathbb{R}$ and $\zeta\in(0,1)$. For $|x-x_0|<\zeta/2$, using $\int_{\mathbb{R}}\Phi_\varepsilon(x-y,\tau)\,dy=1$,
$$
|w(x,\tau)-w_0(x_0)|\le\sup_{|y-x_0|\le\zeta}|w_0(y)-w_0(x_0)|
+\int_{|y-x_0|>\zeta}\Phi_\varepsilon(x-y,\tau)\big(w_0(y)+w_0(x_0)\big)\,dy .
$$
If $|y-x_0|>\zeta$ and $|x-x_0|<\zeta/2$, then $|x-y|\ge|y-x_0|/2$. By Lemma~\ref{lem:U0-growth} with $\eta=1$, $w_0(y)\le e^{(y^2+K_1)/\varepsilon}\le e^{(2(y-x_0)^2+2x_0^2+K_1)/\varepsilon}$. Hence, with $z=y-x_0$ and $0<\tau\le1/32$ (so that $\frac1{8\tau}-2\ge\frac1{16\tau}$),
$$
\int_{|y-x_0|>\zeta}\Phi_\varepsilon(x-y,\tau)\big(w_0(y)+w_0(x_0)\big)dy
\le\frac{C(x_0,\varepsilon)}{\sqrt{2\pi\varepsilon\tau}}\int_{|z|>\zeta}e^{-\frac{z^2}{16\varepsilon\tau}}dz
=\frac{4C(x_0,\varepsilon)}{\sqrt{2\pi}}\int_{|r|>\zeta/(4\sqrt{\varepsilon\tau})}e^{-r^2}dr ,
$$
which tends to $0$ as $\tau\to0^+$. Since $w_0$ is continuous, letting first $(x,\tau)\to(x_0,0^+)$ and then $\zeta\to0^+$ proves the last assertion.
\end{proof}

With the heat equation at hand, the Hopf--Cole transformation produces the viscous solutions. Part (iii) below, which rewrites $v^\varepsilon$ as a weighted average of $u_0$, is the source of the uniform bounds that make the vanishing viscosity limit possible.

\begin{proposition}[Viscous Burgers equation]\label{prop:viscous}
Assume (H3) and let $\varepsilon>0$. Define, for $x\in\mathbb{R}$ and $\tau>0$,
\begin{equation}\label{veps}
v^\varepsilon(x,\tau):=-\varepsilon\,\frac{w_x(x,\tau)}{w(x,\tau)}
=\frac{\displaystyle\int_{\mathbb{R}}\frac{x-y}{\tau}\,e^{-G(x,y,\tau)/\varepsilon}\,dy}
{\displaystyle\int_{\mathbb{R}} e^{-G(x,y,\tau)/\varepsilon}\,dy},
\end{equation}
with $w$ as in Lemma~\ref{lem:heat}. Then the following assertions hold.
\begin{itemize}
\item[(i)] $v^\varepsilon\in C^\infty(\mathbb{R}\times(0,\infty))$ solves \eqref{burgers-visc}.
\item[(ii)] $\lim\limits_{\tau\to0^+}\int_0^xv^\varepsilon(z,\tau)\,dz=U_0(x)$ for every $x\in\mathbb{R}$.
\item[(iii)] If $u_0\in L^\infty(\mathbb{R})$, then
\begin{equation}\label{veps-average}
v^\varepsilon(x,\tau)=\frac{\int_{\mathbb{R}} u_0(y)\,e^{-G(x,y,\tau)/\varepsilon}\,dy}{\int_{\mathbb{R}} e^{-G(x,y,\tau)/\varepsilon}\,dy},
\qquad\text{hence}\qquad |v^\varepsilon(x,\tau)|\le\|u_0\|_{L^\infty}.
\end{equation}
\end{itemize}
\end{proposition}

\begin{proof}
Since $\exp\{-(x-y)^2/(2\varepsilon\tau)\}w_0(y)=e^{-G(x,y,\tau)/\varepsilon}$, differentiating \eqref{heat-sol} under the integral sign gives the second expression in \eqref{veps}.

(i) Let $V:=-\varepsilon\ln w\in C^\infty(\mathbb{R}\times(0,\infty))$. From $w_\tau=\frac\varepsilon2w_{xx}$,
$$
V_\tau=-\frac{\varepsilon^2}{2}\frac{w_{xx}}{w},\qquad V_x=-\varepsilon\frac{w_x}{w},\qquad
V_{xx}=-\varepsilon\Big(\frac{w_{xx}}{w}-\frac{w_x^2}{w^2}\Big),
$$
so that $V_\tau+\frac12V_x^2=\frac\varepsilon2V_{xx}$. Differentiating in $x$ and using $v^\varepsilon=V_x$ gives \eqref{burgers-visc}.

(ii) We have $\int_0^xv^\varepsilon(z,\tau)\,dz=V(x,\tau)-V(0,\tau)$. By Lemma~\ref{lem:heat} and the continuity and positivity of $w_0$, this tends to $-\varepsilon\ln w_0(x)+\varepsilon\ln w_0(0)=U_0(x)-U_0(0)=U_0(x)$.

(iii) Let $F(y):=\exp\{-(x-y)^2/(2\varepsilon\tau)\}$, so that $F'(y)=\frac{x-y}{\varepsilon\tau}F(y)$. Since $U_0$ is Lipschitz, $w_0$ is locally absolutely continuous with $w_0'=-\frac{u_0}{\varepsilon}w_0$ a.e. Integrating by parts on $[-n,n]$,
$$
\int_{-n}^{n}\frac{x-y}{\tau}F(y)w_0(y)\,dy=\varepsilon\big[F\,w_0\big]_{-n}^{n}
+\int_{-n}^{n}F(y)\,u_0(y)\,w_0(y)\,dy .
$$
Since $w_0(y)\le e^{\|u_0\|_{L^\infty}|y|/\varepsilon}$, the boundary term tends to $0$ and both integrals converge absolutely as $n\to\infty$. Dividing by $\int_{\mathbb{R}} Fw_0\,dy$ gives \eqref{veps-average}.
\end{proof}

\section{Minimizers of the variational problem}\label{sec:variational}

As $\varepsilon\to0^+$, the probability measures with densities proportional to $e^{-G(x,\cdot,\tau)/\varepsilon}$ that appear in \eqref{veps} concentrate on the set $\mathcal M(x,\tau)$ of minimizers of $G(x,\cdot,\tau)$. To identify the limit of $v^\varepsilon$ we therefore need precise information on where these minimizers lie, how they depend on $(x,\tau)$, and how often they fail to be unique. The results of this section are the analogues, for initial data satisfying (H3), of the classical properties established by Hopf and Lax \cite{Lax1957}. We include complete proofs because the localization estimates, with explicit constants, are needed to handle unbounded data in the vanishing viscosity limit.

We begin by showing that, for $(x,\tau)$ in a bounded set, all the minimizers lie in a fixed bounded interval.

\begin{lemma}[Localization of minimizers]\label{lem:localization}
Assume (H3), let $R,T>0$, $K:=K_\eta$ with $\eta=1/(16T)$ as in Lemma~\ref{lem:U0-growth}, and
$M_{R,T}:=\max\Big\{2R,\ 4\sqrt{T(\|U_0\|_{L^\infty(-R,R)}+K)}\Big\}$.
Then, for $|x|\le R$, $0<\tau\le T$ and $|y|\ge2R$,
\begin{equation}\label{G-coercive}
G(x,y,\tau)\ge\frac{y^2}{16T}-K,
\end{equation}
and $\mathcal M(x,\tau)\subset[-M_{R,T},M_{R,T}]$.
\end{lemma}

\begin{proof}
Let $|x|\le R$, $0<\tau\le T$ and $|y|\ge2R$. Then $|x-y|\ge|y|-|x|\ge|y|/2$, so $(x-y)^2/(2\tau)\ge y^2/(8T)$, while $U_0(y)\ge-y^2/(16T)-K$. This gives \eqref{G-coercive}. If $|y|>M_{R,T}$, then $|y|>2R$ and $y^2/(16T)>\|U_0\|_{L^\infty(-R,R)}+K$, so $G(x,y,\tau)>\|U_0\|_{L^\infty(-R,R)}\ge U_0(x)=G(x,x,\tau)\ge\min\limits_{z\in\mathbb{R}}G(x,z,\tau)$. Thus $y\notin\mathcal M(x,\tau)$.
\end{proof}

The next lemma contains the key structural property of the minimizers, namely that they move monotonically with $x$. It is the variational counterpart of the fact that characteristics of the Burgers equation do not cross.

\begin{lemma}[Monotonicity]\label{lem:monotone}
Let $\tau>0$, $x_1<x_2$, $y_1\in\mathcal M(x_1,\tau)$ and $y_2\in\mathcal M(x_2,\tau)$. Then $y_1\le y_2$. Consequently $y_-(x_1,\tau)\le y_+(x_1,\tau)\le y_-(x_2,\tau)\le y_+(x_2,\tau)$.
\end{lemma}

\begin{proof}
Adding $G(x_1,y_1,\tau)\le G(x_1,y_2,\tau)$ and $G(x_2,y_2,\tau)\le G(x_2,y_1,\tau)$, the terms involving $U_0$ cancel and we obtain $(x_1-y_1)^2+(x_2-y_2)^2\le(x_1-y_2)^2+(x_2-y_1)^2$, that is, $2(x_2-x_1)(y_1-y_2)\le0$. Since $x_2>x_1$, $y_1\le y_2$.
\end{proof}

We now show that the minimizer is unique outside a null set, which is what makes the limit of $v^\varepsilon$ well defined almost everywhere.

\begin{lemma}[Structure of the minimizers]\label{lem:structure}
Assume (H3).
\begin{itemize}
\item[(i)] $y_+$ is upper semicontinuous and $y_-$ is lower semicontinuous on $\mathbb{R}\times(0,\infty)$. In particular, both are Borel functions.
\item[(ii)] For each $\tau>0$, $x\mapsto y_+(x,\tau)$ is nondecreasing and right-continuous, and $x\mapsto y_-(x,\tau)$ is nondecreasing and left-continuous.
\item[(iii)] The set $N:=\{(x,\tau)\in\mathbb{R}\times(0,\infty):\ y_-(x,\tau)<y_+(x,\tau)\}$, i.e.\ the set where $G(x,\cdot,\tau)$ has more than one minimizer, is Borel, each section $N_\tau:=\{x:(x,\tau)\in N\}$ is countable, and $N$ has Lebesgue measure zero.
\end{itemize}
\end{lemma}

\begin{proof}
(i) Let $(x_n,\tau_n)\to(x,\tau)$ with $\tau>0$, and $y_n:=y_+(x_n,\tau_n)$. For $n$ large, $|x_n|\le|x|+1=:R$ and $\tau_n\le2\tau=:T$, so $|y_n|\le M_{R,T}$ by Lemma~\ref{lem:localization}. Let $y'$ be the limit of a subsequence. Passing to the limit in $G(x_n,y_n,\tau_n)\le G(x_n,z,\tau_n)$, $z\in\mathbb{R}$, and using the continuity of $G$ we get $G(x,y',\tau)\le G(x,z,\tau)$ for all $z$, i.e.\ $y'\in\mathcal M(x,\tau)$, so $y'\le y_+(x,\tau)$. Hence $\limsup\limits_{n\to\infty}y_+(x_n,\tau_n)\le y_+(x,\tau)$. The argument for $y_-$ is symmetric.

(ii) Monotonicity follows from Lemma~\ref{lem:monotone}. If $\bar x\to x^+$, then $y_+(\bar x,\tau)\ge y_+(x,\tau)$ by monotonicity and $\limsup\limits_{\bar x\to x^+}y_+(\bar x,\tau)\le y_+(x,\tau)$ by (i). The statement for $y_-$ is analogous.

(iii) Since $y_+-y_-$ is upper semicontinuous, each set $\{y_+-y_-\ge1/n\}$ is relatively closed in $\mathbb{R}\times(0,\infty)$, and $N=\bigcup\limits_{n\ge1}\{y_+-y_-\ge1/n\}$ is Borel. For fixed $\tau$, the open intervals $(y_-(x,\tau),y_+(x,\tau))$, $x\in N_\tau$, are nonempty and, by Lemma~\ref{lem:monotone}, pairwise disjoint. Choosing a rational number in each of them shows that $N_\tau$ is countable. By Tonelli's theorem, $|N|=\int_0^\infty|N_\tau|\,d\tau=0$.
\end{proof}

The quantity $(x-y)/\tau$, with $y$ a minimizer, is the velocity of the characteristic that reaches $(x,\tau)$. The following bounds on it are the a priori estimates of the problem. They give local boundedness for general data, global boundedness for bounded data, and decay for integrable data.

\begin{lemma}[Bounds on the velocity]\label{lem:velocity}
Assume (H3), let $x\in\mathbb{R}$, $\tau>0$ and $y\in\mathcal M(x,\tau)$.
\begin{itemize}
\item[(i)] If $|x|\le R$ and $\tau\le T$, then $|x-y|\le L_{R,T}\,\tau$, where $L_{R,T}:=\|u_0\|_{L^\infty(-M_{R,T}-1,\,M_{R,T}+1)}$.
\item[(ii)] If $u_0\in L^\infty(\mathbb{R})$, then $|x-y|\le\|u_0\|_{L^\infty}\,\tau$.
\item[(iii)] If $u_0\in L^1(\mathbb{R})$, then $(x-y)^2\le2\tau\|u_0\|_{L^1}$.
\end{itemize}
\end{lemma}

\begin{proof}
(i) By Lemma~\ref{lem:localization}, $y\in[-M_{R,T},M_{R,T}]$. For $0<|z|<1$, the inequality $G(x,y,\tau)\le G(x,y+z,\tau)$ reads $2z(x-y)-z^2\le2\tau\big(U_0(y+z)-U_0(y)\big)\le2\tau L_{R,T}|z|$. Dividing by $2|z|$ we get $(x-y)-\frac z2\le\tau L_{R,T}$ if $z>0$, and $-(x-y)-\frac{|z|}2\le\tau L_{R,T}$ if $z<0$. Letting $z\to0$ gives the claim.

(ii) The same argument applies with $L_{R,T}$ replaced by the global Lipschitz constant $\|u_0\|_{L^\infty}$ of $U_0$.

(iii) From $G(x,y,\tau)\le G(x,x,\tau)$ we get $(x-y)^2/(2\tau)\le U_0(x)-U_0(y)=\int_y^xu_0(r)\,dr\le\|u_0\|_{L^1}$.
\end{proof}

For unbounded data, the viscous solutions need not be bounded uniformly in $\varepsilon$, and the uniform bound \eqref{veps-average} is not available. The following lemma circumvents this difficulty. On a bounded set of $(x,\tau)$, the variational problem only sees the initial datum on a bounded interval, so that $u_0$ may be replaced by a bounded and integrable function without changing the minimizers there.

\begin{lemma}\label{lem:cutoff}
Assume (H3), let $R,T>0$, $M:=M_{R,T}$, and $\tilde u_0:=u_0\mathbf 1_{[-M,M]}$. Denote by $\tilde U_0$, $\tilde G$, $\tilde{\mathcal M}$ the objects \eqref{G} associated with $\tilde u_0$. Then $\tilde u_0\in L^\infty(\mathbb{R})\cap L^1(\mathbb{R})$ satisfies (H3), $\tilde u_0=u_0$ on $[-R,R]$, and $\tilde{\mathcal M}(x,\tau)=\mathcal M(x,\tau)$ for $|x|\le R$ and $0<\tau\le T$.
\end{lemma}

\begin{proof}
The first assertions follow from $M\ge2R$ and the boundedness of $\tilde U_0$. Fix $|x|\le R$ and $0<\tau\le T$, and set $\underline G:=\min\limits_{y\in\mathbb{R}}G(x,y,\tau)$. On $[-M,M]$ we have $\tilde G=G\ge\underline G$. If $y>M$, then $\tilde U_0(y)=U_0(M)$ and, since $x\le R<M<y$,
$$
\tilde G(x,y,\tau)=\frac{(x-y)^2}{2\tau}+U_0(M)>\frac{(x-M)^2}{2\tau}+U_0(M)=G(x,M,\tau)\ge\underline G.
$$
The case $y<-M$ is symmetric. Hence $\tilde G\ge\underline G$ on $\mathbb{R}$, with equality exactly at the points of $\mathcal M(x,\tau)$, which lie in $[-M,M]$ by Lemma~\ref{lem:localization}. Therefore $\tilde{\mathcal M}(x,\tau)=\mathcal M(x,\tau)$.
\end{proof}

\section{Vanishing viscosity limit}\label{sec:limit}

We now let $\varepsilon\to0^+$. The argument has two steps. At points where the minimizer is unique, a Laplace-type argument shows that the weighted average \eqref{veps} converges to the velocity of the unique minimizing characteristic, and by Lemma~\ref{lem:structure} this happens almost everywhere. The limit is then shown to be a weak solution by passing to the limit in the weak formulation of the viscous problem, after a localization that reduces the problem to bounded data. Finally, Lemmas~\ref{lem:weak-reduction} and~\ref{lem:viscous-reduction} carry these results over to \eqref{main} through the change of variables \eqref{transf}.

\begin{lemma}\label{lem:laplace}
Assume (H3), and let $(x,\tau)\in(\mathbb{R}\times(0,\infty))\setminus N$ with unique minimizer $y^*$. Then $\lim\limits_{\varepsilon\to0^+}v^\varepsilon(x,\tau)=(x-y^*)/\tau$.
\end{lemma}

\begin{proof}
Let $P(y):=G(x,y,\tau)-G(x,y^*,\tau)$. The function $P:\mathbb{R}\to[0,\infty)$ is continuous and vanishes only at $y=y^*$. By \eqref{G-coercive}, $P$ grows quadratically, hence $I:=\int_{\mathbb{R}}|y-y^*|e^{-P(y)}dy<\infty$. Multiplying numerator and denominator in \eqref{veps} by $e^{G(x,y^*,\tau)/\varepsilon}$,
$$
\Big|v^\varepsilon(x,\tau)-\frac{x-y^*}{\tau}\Big|
=\frac{\big|\int_{\mathbb{R}}\frac{y^*-y}{\tau}e^{-P/\varepsilon}dy\big|}{\int_{\mathbb{R}} e^{-P/\varepsilon}dy}
\le\frac1\tau\,\frac{\int_{\mathbb{R}}|y-y^*|e^{-P/\varepsilon}dy}{\int_{\mathbb{R}} e^{-P/\varepsilon}dy}.
$$
Fix $\zeta>0$. Since $P$ is continuous, positive on the closed set $\{|y-y^*|\ge\zeta\}$ and tends to $+\infty$ as $|y|\to\infty$, we have $\gamma:=\frac12\inf\{P(y):|y-y^*|\ge\zeta\}>0$. Choose $\zeta'\in(0,\zeta)$ with $P\le\gamma$ on $[y^*-\zeta',y^*+\zeta']$. For $0<\varepsilon\le1$ and $|y-y^*|\ge\zeta$ we have $P\ge2\gamma$, hence $P/\varepsilon\ge(P-2\gamma)+2\gamma/\varepsilon$ and $e^{-P/\varepsilon}\le e^{2\gamma}e^{-2\gamma/\varepsilon}e^{-P}$. Therefore
$$
\int_{\mathbb{R}}|y-y^*|e^{-P/\varepsilon}dy\le\zeta\int_{\mathbb{R}} e^{-P/\varepsilon}dy+I\,e^{2\gamma}e^{-2\gamma/\varepsilon},
\qquad
\int_{\mathbb{R}} e^{-P/\varepsilon}dy\ge2\zeta'e^{-\gamma/\varepsilon},
$$
and $\big|v^\varepsilon(x,\tau)-\frac{x-y^*}{\tau}\big|\le\frac\zeta\tau+\frac{I\,e^{2\gamma}}{2\zeta'\tau}e^{-\gamma/\varepsilon}$. Letting $\varepsilon\to0^+$ and then $\zeta\to0^+$ gives the claim.
\end{proof}

With the pointwise limit identified, we can carry out the vanishing viscosity limit for the classical Burgers equation with initial data satisfying (H3).

\begin{theorem}[Vanishing viscosity for the Burgers equation]\label{thm:burgers}
Assume (H3) and define $v(x,\tau):=(x-y_+(x,\tau))/\tau$ for $x\in\mathbb{R}$, $\tau>0$. Then the following assertions hold.
\begin{itemize}
\item[(i)] $v$ is Borel and $v^\varepsilon\to v$ on $(\mathbb{R}\times(0,\infty))\setminus N$ as $\varepsilon\to0^+$.
\item[(ii)] $|v|\le L_{R,T}$ on $[-R,R]\times(0,T]$ for all $R,T>0$. If $u_0\in L^\infty(\mathbb{R})$, then $|v|\le\|u_0\|_{L^\infty}$, and if $u_0\in L^1(\mathbb{R})$, then $|v(x,\tau)|\le(2\|u_0\|_{L^1}/\tau)^{1/2}$.
\item[(iii)] $v(x_2,\tau)-v(x_1,\tau)\le(x_2-x_1)/\tau$ for $\tau>0$ and $x_1<x_2$. Moreover, $v(\cdot,\tau)\in BV_{\mathrm{loc}}(\mathbb{R})$ is right-continuous and $v(x-,\tau)>v(x+,\tau)$ at each of its countably many discontinuity points.
\item[(iv)] $v\in L^\infty_{\mathrm{loc}}(\mathbb{R}\times[0,\infty))$ is a weak solution of \eqref{burgers} in $\mathbb{R}\times(0,\infty)$, and hence in $\mathbb{R}\times(0,T^*)$ for every $T^*\in(0,\infty]$.
\end{itemize}
\end{theorem}

\begin{proof}
(i) follows from Lemmas~\ref{lem:structure} and~\ref{lem:laplace}, and (ii) from Lemma~\ref{lem:velocity}. (iii) Since $y_+(\cdot,\tau)$ is nondecreasing, $v(x_2,\tau)-v(x_1,\tau)=\frac{x_2-x_1}{\tau}-\frac{y_+(x_2,\tau)-y_+(x_1,\tau)}{\tau}\le\frac{x_2-x_1}{\tau}$. Moreover $v(\cdot,\tau)$ is the difference of the smooth function $x/\tau$ and the nondecreasing right-continuous function $y_+(\cdot,\tau)/\tau$, which yields the remaining assertions.

(iv) Local boundedness follows from (ii). Fix $\psi\in C_c^\infty(\mathbb{R}\times[0,\infty))$ and choose $R,T>0$ with $\operatorname{supp}\psi\subset(-R,R)\times[0,T)$. Let $\tilde u_0$ be the cut-off datum of Lemma~\ref{lem:cutoff}, let $\tilde v^\varepsilon$ be the corresponding viscous solution \eqref{veps}, and $\tilde v:=(x-\tilde y_+)/\tau$ the corresponding limit. By Lemma~\ref{lem:cutoff}, $\tilde v=v$ on $[-R,R]\times(0,T]$, and $\tilde u_0=u_0$ on $[-R,R]$. By Proposition~\ref{prop:viscous}(iii), $|\tilde v^\varepsilon|\le\tilde C:=\|\tilde u_0\|_{L^\infty}$ for all $\varepsilon>0$.

Since $\tilde v^\varepsilon$ is a smooth solution of \eqref{burgers-visc}, multiplying by $\psi$, integrating over $\mathbb{R}\times(\tau_1,\infty)$, $\tau_1\in(0,T)$, and integrating by parts we get
\begin{equation}\label{visc-weak-tau1}
\int_{\tau_1}^\infty\!\!\int_{\mathbb{R}}\Big[\tilde v^\varepsilon\psi_\tau+\frac{(\tilde v^\varepsilon)^2}{2}\psi_x
+\frac\varepsilon2\tilde v^\varepsilon\psi_{xx}\Big]dx\,d\tau+\int_{\mathbb{R}}\tilde v^\varepsilon(x,\tau_1)\psi(x,\tau_1)\,dx=0.
\end{equation}
Let $\tilde V^\varepsilon(x,\tau):=\int_0^x\tilde v^\varepsilon(z,\tau)\,dz$. Then
$$
\int_{\mathbb{R}}\tilde v^\varepsilon(x,\tau_1)\psi(x,\tau_1)\,dx=-\int_{\mathbb{R}}\tilde V^\varepsilon(x,\tau_1)\psi_x(x,\tau_1)\,dx,
$$
$|\tilde V^\varepsilon(x,\tau_1)|\le\tilde C|x|$, $\tilde V^\varepsilon(x,\tau_1)\to\tilde U_0(x)$ for every $x$ by Proposition~\ref{prop:viscous}(ii), and $\psi_x(\cdot,\tau_1)\to\psi_x(\cdot,0)$ uniformly. By dominated convergence on $[-R,R]$ and an integration by parts ($\tilde U_0$ is Lipschitz),
$$
\lim_{\tau_1\to0^+}\int_{\mathbb{R}}\tilde v^\varepsilon(x,\tau_1)\psi(x,\tau_1)\,dx
=-\int_{\mathbb{R}}\tilde U_0(x)\psi_x(x,0)\,dx=\int_{\mathbb{R}}\tilde u_0(x)\psi(x,0)\,dx .
$$
The double integral in \eqref{visc-weak-tau1} converges as $\tau_1\to0^+$ because its integrand is bounded and compactly supported. Hence
\begin{equation}\label{visc-weak}
\int_0^\infty\!\!\int_{\mathbb{R}}\Big[\tilde v^\varepsilon\psi_\tau+\frac{(\tilde v^\varepsilon)^2}{2}\psi_x\Big]dx\,d\tau
+\int_{\mathbb{R}}\tilde u_0\psi(x,0)\,dx=-\frac\varepsilon2\int_0^\infty\!\!\int_{\mathbb{R}}\tilde v^\varepsilon\psi_{xx}\,dx\,d\tau .
\end{equation}
By (i) applied to $\tilde u_0$, $\tilde v^\varepsilon\to\tilde v$ a.e.\ in $\mathbb{R}\times(0,\infty)$, and $|\tilde v^\varepsilon|\le\tilde C$. By dominated convergence, letting $\varepsilon\to0^+$ in \eqref{visc-weak} yields $\int_0^\infty\!\int_{\mathbb{R}}[\tilde v\psi_\tau+\frac12\tilde v^2\psi_x]\,dx\,d\tau+\int_{\mathbb{R}}\tilde u_0\psi(x,0)\,dx=0$. Since $\tilde v=v$ and $\tilde u_0=u_0$ on the support of $\psi$, \eqref{weak} holds for $v$ (with $\alpha\equiv1$, $\sigma\equiv0$). Finally, every $\psi\in C_c^\infty(\mathbb{R}\times[0,T^*))$ extends by zero to an element of $C_c^\infty(\mathbb{R}\times[0,\infty))$.
\end{proof}

It remains to deduce Theorem~\ref{thm:existence} from Theorem~\ref{thm:burgers} by means of the change of variables \eqref{transf}. All the properties of $v$ are preserved by the multiplication by $e^{-S(t)}$ and the time change $\tau=A(t)$. The only point requiring attention is that the null set $N$ is mapped onto a null set.

\begin{proof}[Proof of Theorem~\ref{thm:existence}]
(i) Let $v^\varepsilon$ be given by \eqref{veps}. Then $u^\varepsilon(x,t)=e^{-S(t)}v^\varepsilon(x,A(t))$ is the function \eqref{ueps}. Its regularity and the fact that it solves \eqref{viscous} follow from Proposition~\ref{prop:viscous}(i) and Lemma~\ref{lem:viscous-reduction}. Moreover $\int_0^xu^\varepsilon(z,t)\,dz=e^{-S(t)}\int_0^xv^\varepsilon(z,A(t))\,dz\to U_0(x)$ as $t\to0^+$, by Proposition~\ref{prop:viscous}(ii), since $S(t)\to0$ and $A(t)\to0^+$.

(ii) We have $u(x,t)=e^{-S(t)}v(x,A(t))$ with $v$ as in Theorem~\ref{thm:burgers}. Hence $u$ is Borel. Let $\mathcal N:=\Lambda^{-1}\big(N\cap(\mathbb{R}\times(0,A_\infty))\big)$ with $\Lambda(x,t)=(x,A(t))$. Since $\Lambda$ is a $C^1$ diffeomorphism, $\mathcal N$ is a Borel null set, and $\mathcal N_t=N_{A(t)}$ is countable by Lemma~\ref{lem:structure}. For $(x,t)\notin\mathcal N$, Theorem~\ref{thm:burgers}(i) gives $u^\varepsilon(x,t)=e^{-S(t)}v^\varepsilon(x,A(t))\to e^{-S(t)}v(x,A(t))=u(x,t)$.

(iii) For $|x|\le R$ and $0<t\le T$ we have $A(t)\le A(T)$, and Theorem~\ref{thm:burgers}(ii) gives $|u(x,t)|\le e^{-S(t)}L_{R,A(T)}$.

(iv) By Theorem~\ref{thm:burgers}(iv), $v$ is a weak solution of \eqref{burgers} in $\mathbb{R}\times(0,A_\infty)$, and $u$ is related to $v$ by \eqref{v-from-u}. The claim follows from Lemma~\ref{lem:weak-reduction}.

(v) This follows from Theorem~\ref{thm:burgers}(iii) with $\tau=A(t)$, multiplied by $e^{-S(t)}>0$.
\end{proof}

\section{Uniqueness and decay of entropy solutions}\label{sec:uniqueness}

In this section we prove Theorem~\ref{thm:uniqueness}. The decay estimates follow at once from the velocity bounds of Lemma~\ref{lem:velocity}. The heart of the matter is uniqueness. By Lemmas~\ref{lem:weak-reduction} and~\ref{lem:onesided}, the change of variables \eqref{transf} maps entropy solutions of \eqref{main} in $\mathbb{R}\times(0,T)$ onto entropy solutions of the classical Burgers equation \eqref{burgers} in $\mathbb{R}\times(0,A(T))$. Hence it suffices to prove uniqueness for \eqref{burgers} on a strip of finite height.

We use the duality argument of Oleinik \cite{Oleinik1957}, in the form presented in \cite[Section~3.4.3]{Evans2010}. Its idea is as follows. If $v_1$ and $v_2$ are two weak solutions with the same initial datum, then $W:=v_1-v_2$ is a weak solution of the linear equation $W_\tau+(bW)_x=0$ with $b:=\frac12(v_1+v_2)$ and zero initial datum. To show that $W=0$, one tests this equation with a solution $\varphi$ of the adjoint backward transport problem $\varphi_\tau+b\varphi_x=\psi$, $\varphi(\cdot,T)=0$, with $\psi$ arbitrary, which would give $\iint W\psi=0$. Since $b$ is discontinuous, the adjoint problem is solved with a regularized coefficient $b_\delta$, and an error term $\iint W(b-b_\delta)\varphi_x$ has to be controlled. This is where the entropy condition enters. An upper bound on $(b_\delta)_x$ means that backward characteristics of the adjoint problem do not concentrate, which yields a bound on $\varphi_x$ that is uniform in $\delta$. Without the entropy condition, $\varphi_x$ may concentrate on the discontinuities of $b$ and the argument fails, in agreement with the nonuniqueness of general weak solutions.

\begin{proposition}[Oleinik]\label{prop:oleinik}
Let $u_0\in L^\infty(\mathbb{R})$ and $T^*\in(0,\infty)$. Let $v_1,v_2\in L^\infty(\mathbb{R}\times(0,T^*))$ be weak solutions of \eqref{burgers} in $\mathbb{R}\times(0,T^*)$ such that, for some $C\ge0$,
\begin{equation}\label{oleinik-vi}
v_i(x+z,\tau)-v_i(x,\tau)\le C\Big(1+\frac1\tau\Big)z\quad
\text{for a.e. }(x,z,\tau)\in\mathbb{R}\times(0,\infty)\times(0,T^*),\ i=1,2.
\end{equation}
Then $v_1=v_2$ a.e.\ in $\mathbb{R}\times(0,T^*)$.
\end{proposition}

\begin{proof}
Let $B\ge\max\{\|v_1\|_{L^\infty},\|v_2\|_{L^\infty}\}$, $W:=v_1-v_2$ and $b:=\frac12(v_1+v_2)$, so that $|W|\le2B$, $|b|\le B$ and $\frac12v_1^2-\frac12v_2^2=bW$.

We first derive a linear equation for $W$. By Lemma~\ref{lem:lipschitz}, the weak formulation holds for $v_1$ and $v_2$ with every Lipschitz test function $\varphi$ with compact support in $\mathbb{R}\times[0,T^*)$. Subtracting the two identities, the terms containing $u_0$ cancel, and
\begin{equation}\label{dual-identity}
\int_0^{T^*}\!\!\int_{\mathbb{R}} W\big(\varphi_\tau+b\,\varphi_x\big)\,dx\,d\tau=0
\end{equation}
for every such $\varphi$.

We next regularize the coefficient $b$. Fix $\psi\in C_c^\infty(\mathbb{R}\times(0,T^*))$, and choose $R>0$ and $T\in(0,T^*)$ with $\operatorname{supp}\psi\subset[-R,R]\times(0,T)$. We shall prove that $\iint W\psi=0$, which gives $W=0$ a.e.\ because $\psi$ is arbitrary.

Extend $b$ by zero outside $\mathbb{R}\times(0,T^*)$, let $J_\delta$ be a standard mollifier on $\mathbb{R}^2$ supported in the ball of radius $\delta$, and set
$$
b_\delta:=J_\delta*b,\qquad 0<\delta<\delta_0:=\min\{T^*-T,\,T/2\}.
$$
Then $b_\delta\in C^\infty(\mathbb{R}^2)$, $|b_\delta|\le B$, all derivatives of $b_\delta$ are bounded on $\mathbb{R}^2$ (by constants depending on $\delta$), and $b_\delta\to b$ a.e.\ in $\mathbb{R}\times(0,T^*)$ as $\delta\to0^+$, since $b_\delta\to b$ at every Lebesgue point of $b$.

We claim that
\begin{equation}\label{bdelta-onesided}
(b_\delta)_x(x,\tau)\le2C\Big(1+\frac1\tau\Big)\qquad\text{for all }x\in\mathbb{R},\ 2\delta\le\tau\le T.
\end{equation}
Since $b$ is the average of $v_1$ and $v_2$, it satisfies \eqref{oleinik-vi} as well. By Fubini's theorem, there is a null set $Z\subset(0,\infty)$ such that, for every $z\in(0,\infty)\setminus Z$, the inequality $b(x+z,\tau)-b(x,\tau)\le C(1+1/\tau)z$ holds for a.e.\ $(x,\tau)\in\mathbb{R}\times(0,T^*)$. Fix such a $z$ and $(x,\tau)$ with $2\delta\le\tau\le T$. For $|(y,r)|<\delta$ we have $\tau-r\in(\tau/2,T^*)$, hence $1/(\tau-r)\le2/\tau$, and
$$
b_\delta(x+z,\tau)-b_\delta(x,\tau)
=\int\!\!\int J_\delta(y,r)\big[b(x+z-y,\tau-r)-b(x-y,\tau-r)\big]\,dy\,dr
\le C\Big(1+\frac2\tau\Big)z .
$$
Since $b_\delta$ is continuous and $(0,\infty)\setminus Z$ is dense, this inequality holds for every $z>0$. Dividing by $z$ and letting $z\to0^+$ gives $(b_\delta)_x\le C(1+2/\tau)\le2C(1+1/\tau)$, which is \eqref{bdelta-onesided}.

We now solve the adjoint problem with the regularized coefficient. For $(x,\tau)\in\mathbb{R}\times[0,T]$, let $r\mapsto X(r;x,\tau)$, $r\in[0,T]$, be the solution of
$
\frac{d}{dr}X(r;x,\tau)=b_\delta\big(X(r;x,\tau),r\big)$ with $X(\tau;x,\tau)=x$.
Since $b_\delta$ is smooth, bounded and globally Lipschitz in $x$, $X$ exists on $[0,T]$, is of class $C^\infty$ in $(r,x,\tau)$, satisfies $|X(r;x,\tau)-x|\le B|r-\tau|$, and has the flow property $X(r';X(r;x,\tau),r)=X(r';x,\tau)$. Define
$
\varphi^\delta(x,\tau):=-\int_\tau^T\psi\big(X(r;x,\tau),r\big)\,dr$ for $0\le\tau\le T$ and
$\varphi^\delta(x,\tau):=0$ for $T<\tau<T^*$.
Then $\varphi^\delta\in C^\infty(\mathbb{R}\times[0,T])$. By the flow property, $\varphi^\delta(X(r;x,\tau),r)=-\int_r^T\psi(X(r';x,\tau),r')\,dr'$. Differentiating in $r$ at $r=\tau$ gives
\begin{equation}\label{adjoint}
\varphi^\delta_\tau+b_\delta\,\varphi^\delta_x=\psi\quad\text{in }\mathbb{R}\times[0,T],
\qquad \varphi^\delta(\cdot,T)=0 .
\end{equation}
If $|x|>R':=R+BT$, then $|X(r;x,\tau)|>R$ for all $r\in[0,T]$, so $\varphi^\delta(x,\tau)=0$. Hence $\varphi^\delta$ is supported in $[-R',R']\times[0,T]$, a compact set independent of $\delta$. Since $\varphi^\delta$ is smooth on $\mathbb{R}\times[0,T]$, vanishes at $\tau=T$ and vanishes for $\tau>T$, it is Lipschitz continuous with compact support in $\mathbb{R}\times[0,T^*)$, and it is an admissible test function in \eqref{dual-identity}. Using \eqref{adjoint} we obtain
\begin{equation}\label{error}
\int_0^{T^*}\!\!\int_{\mathbb{R}} W\psi\,dx\,d\tau
=-\int_0^{T}\!\!\int_{\mathbb{R}} W\,(b-b_\delta)\,\varphi^\delta_x\,dx\,d\tau .
\end{equation}

The right-hand side of \eqref{error} is controlled by two estimates for $q:=\varphi^\delta_x$. Differentiating \eqref{adjoint} in $x$, we obtain
\begin{equation}\label{q}
q_\tau+b_\delta q_x+(b_\delta)_x\,q=\psi_x\quad\text{in }\mathbb{R}\times[0,T],\qquad q(\cdot,T)=0 .
\end{equation}

The first estimate is an $L^\infty$ bound away from $\tau=0$. Fix $(x,\tau)$ with $2\delta\le\tau\le T$, and let $\bar q(r):=q(X(r;x,\tau),r)$ for $r\in[\tau,T]$. By \eqref{q}, $\bar q'(r)=-(b_\delta)_x(X(r;x,\tau),r)\,\bar q(r)+\psi_x(X(r;x,\tau),r)$ and $\bar q(T)=0$, whence
$$
q(x,\tau)=\bar q(\tau)=-\int_\tau^T\psi_x\big(X(r;x,\tau),r\big)\exp\Big(\int_\tau^r(b_\delta)_x\big(X(r';x,\tau),r'\big)\,dr'\Big)dr .
$$
By \eqref{bdelta-onesided}, $(b_\delta)_x(X(r';x,\tau),r')\le2C(1+1/r')$ for $r'\ge\tau\ge2\delta$. Hence, for $\tau_0\in[2\delta,T]$ and $\tau\in[\tau_0,T]$,
\begin{equation}\label{q-Linf}
|q(x,\tau)|\le\|\psi_x\|_{L^\infty}\int_\tau^Te^{2C(r-\tau)}\Big(\frac r\tau\Big)^{2C}dr
\le D_{\tau_0}:=T\|\psi_x\|_{L^\infty}e^{2CT}\Big(\frac T{\tau_0}\Big)^{2C}.
\end{equation}
Only the upper bound on $(b_\delta)_x$ has been used, and the constant $D_{\tau_0}$ does not depend on $\delta$.

The second estimate is an $L^1$ bound up to $\tau=0$. Equation \eqref{q} can be written in conservation form, $q_\tau+(b_\delta q)_x=\psi_x$. For $\omega>0$ let $H_\omega(q):=\sqrt{q^2+\omega^2}-\omega$, so that $0\le H_\omega(q)\le|q|$ and $H_\omega(q)$ increases to $|q|$ as $\omega\to0^+$. Since $q(\cdot,\tau)$ is smooth and supported in $[-R',R']$ for every $\tau\in[0,T]$,
$$
\frac{d}{d\tau}\int_{\mathbb{R}} H_\omega(q)\,dx=\int_{\mathbb{R}} H_\omega'(q)\psi_x\,dx
-\int_{\mathbb{R}} H_\omega'(q)(b_\delta q)_x\,dx .
$$
Integrating by parts, $\int_{\mathbb{R}} H_\omega'(q)(b_\delta q)_x\,dx=\int_{\mathbb{R}}(b_\delta)_x\big[H_\omega'(q)q-H_\omega(q)\big]dx$, and $0\le H_\omega'(q)q-H_\omega(q)=\omega-\omega^2/\sqrt{q^2+\omega^2}\le\omega$, with equality to zero where $q=0$. Since $|H_\omega'|\le1$ and $(b_\delta)_x$ is bounded,
$$
\frac{d}{d\tau}\int_{\mathbb{R}} H_\omega(q)\,dx\ge-\|\psi_x(\cdot,\tau)\|_{L^1}-2R'\|(b_\delta)_x\|_{L^\infty}\,\omega .
$$
Integrating from $\tau$ to $T$ and using $q(\cdot,T)=0$, we get $\int_{\mathbb{R}} H_\omega(q(x,\tau))\,dx\le\int_\tau^T\|\psi_x(\cdot,r)\|_{L^1}dr+2R'T\|(b_\delta)_x\|_{L^\infty}\,\omega$. Letting $\omega\to0^+$, by monotone convergence,
\begin{equation}\label{q-L1}
\int_{\mathbb{R}}|q(x,\tau)|\,dx\le\|\psi_x\|_{L^1(\mathbb{R}\times(0,T))}\qquad\text{for all }\tau\in[0,T].
\end{equation}
This bound holds without any one-sided condition, but it does not prevent $q$ from concentrating, which is why \eqref{q-Linf} is needed.

We can now conclude. Fix $\tau_0\in(0,T)$ and let $\delta<\min\{\delta_0,\tau_0/2\}$. We split the right-hand side of \eqref{error} at $\tau=\tau_0$. Using $|W|\le2B$, $|b-b_\delta|\le2B$ and \eqref{q-L1},
$$
\int_0^{\tau_0}\!\!\int_{\mathbb{R}}|W|\,|b-b_\delta|\,|q|\,dx\,d\tau\le4B^2\|\psi_x\|_{L^1(\mathbb{R}\times(0,T))}\,\tau_0 .
$$
Using \eqref{q-Linf} and the fact that $q$ vanishes for $|x|>R'$,
$$
\int_{\tau_0}^T\!\!\int_{\mathbb{R}}|W|\,|b-b_\delta|\,|q|\,dx\,d\tau\le2BD_{\tau_0}\int_{\tau_0}^T\!\!\int_{-R'}^{R'}|b-b_\delta|\,dx\,d\tau
\longrightarrow0\qquad\text{as }\delta\to0^+,
$$
by dominated convergence, since $b_\delta\to b$ a.e.\ and $|b-b_\delta|\le2B$. Therefore $$ \left|\iint W\psi \right|\le4B^2 \tau_0 \|\psi_x\|_{L^1(\mathbb{R}\times(0,T))}$$ for every $\tau_0\in(0,T)$, and letting $\tau_0\to0^+$ gives $\iint W\psi=0$.
\end{proof}

The proof uses the entropy condition only through the upper bound \eqref{bdelta-onesided}, and only for $\tau\ge2\delta$. Near $\tau=0$ the weaker bound \eqref{q-L1}, valid for any bounded coefficient, suffices because the strip $(0,\tau_0)$ is thin. The singular factor $1/\tau$ in \eqref{oleinik-vi} is therefore compatible with the argument, which is essential since rarefaction waves satisfy the one-sided bound only with this factor. We can now prove Theorem~\ref{thm:uniqueness}.

\begin{proof}[Proof of Theorem~\ref{thm:uniqueness}]
If $u_0\in L^\infty(\mathbb{R})$, then $|U_0(x)|\le\|u_0\|_{L^\infty}|x|=o(x^2)$, so (H3) holds. Assertions (i) and (iii) follow from Theorem~\ref{thm:burgers}(ii) with $\tau=A(t)$, after multiplication by $e^{-S(t)}$.

(ii) By Theorem~\ref{thm:existence}, $u$ is a weak solution, and by (i) it is bounded on $\mathbb{R}\times(0,T)$ for every $T>0$. Moreover, for $0<t<T$, \eqref{oleinik-u} and $A(t)\ge\underline a_Tt$ give $u(x+z,t)-u(x,t)\le e^{\bar S_T}\underline a_T^{-1}z/t$. Hence $u$ is an entropy solution.
Let $\hat u$ be another entropy solution and fix $T>0$. Define $\hat v$ from $\hat u$, and $v$ from $u$, by \eqref{v-from-u}. By Lemma~\ref{lem:weak-reduction}, $\hat v$ and $v$ are weak solutions of \eqref{burgers} in $\mathbb{R}\times(0,A_\infty)$, hence in $\mathbb{R}\times(0,A(T))$, since every element of $C_c^\infty(\mathbb{R}\times[0,A(T)))$ belongs to $C_c^\infty(\mathbb{R}\times[0,A_\infty))$. Both are bounded in $\mathbb{R}\times(0,A(T))$, because $e^{S}$ is bounded on $[0,T]$, and by Lemma~\ref{lem:onesided}(ii) both satisfy \eqref{oleinik-vi} in $\mathbb{R}\times(0,A(T))$ with a common constant. Proposition~\ref{prop:oleinik} with $T^*=A(T)$ yields $\hat v=v$ a.e.\ in $\mathbb{R}\times(0,A(T))$, that is, $\hat u=u$ a.e.\ in $\mathbb{R}\times(0,T)$. Since $T$ is arbitrary, $\hat u=u$ a.e.\ in $\mathbb{R}\times(0,\infty)$.
\end{proof}

\section{Shock formation from smooth initial data}\label{sec:threshold}

In this section we answer the question, raised in the Introduction, of whether damping and the time dependence of the coefficient $\alpha$ can prevent the formation of shocks from smooth initial data. A shock forms from smooth data when the solution, while remaining bounded, develops an infinitely steep front, that is, $\sup\limits_{x\in\mathbb{R}}|u_x(x,t)|\to\infty$ as $t$ increases to some finite time $t^*$. In acoustics this is the breaking of the waveform \cite{Whitham1974}. For the classical Burgers equation $v_\tau+vv_x=0$, characteristics are the straight lines $x=\xi+\tau u_0(\xi)$, and two neighboring characteristics issued from a region where $u_0$ has slope $-\kappa<0$ meet at $\tau=1/\kappa$. Hence a shock forms exactly at $\tau=1/\kappa$, where $-\kappa$ is the most negative slope of $u_0$ \cite{Lax1964}.

For \eqref{main}, the characteristics are given by \eqref{char-intro}. They are the characteristics of the Burgers equation with $\tau$ replaced by $A(t)$, while the damping only multiplies the solution by the positive factor $e^{-S(t)}$ and does not move them. Two characteristics issued from a region where $u_0$ has slope $-\kappa$ therefore meet when $A(t)$ reaches $1/\kappa$. Since $A(t)$ increases to the maximal displacement $A_\infty$, this happens in finite time if and only if $1/\kappa<A_\infty$. In other words, a shock forms if and only if the characteristics can travel far enough to cross. We now make this precise, and then illustrate it with examples in which the outcome is determined by the damping, by the coefficient $\alpha$, or by their interplay.

\begin{proof}[Proof of Theorem~\ref{thm:threshold}]
Let $\tau^*:=1/\kappa$ if $\kappa>0$ and $\tau^*:=\infty$ if $\kappa=0$. Since $u_0\in C^1(\mathbb{R})$, $U_0\in C^2(\mathbb{R})$ and, for $0<\tau<\tau^*$,
$$
G_{yy}(x,y,\tau)=\frac1\tau+u_0'(y)\ge\frac1\tau-\kappa>0 .
$$
Thus $G(x,\cdot,\tau)$ is strictly convex and coercive, and its unique minimizer $y$ is characterized by $G_y(x,y,\tau)=0$, i.e.\ by
\begin{equation}\label{X}
x=\chi(y,\tau):=y+\tau u_0(y).
\end{equation}
For $0\le\tau<\tau^*$, $\chi_y=1+\tau u_0'\ge1-\kappa\tau>0$ and $|\chi(y,\tau)-y|\le\tau\|u_0\|_{L^\infty}$, so $\chi(\cdot,\tau):\mathbb{R}\to\mathbb{R}$ is an increasing bijection. Let $\xi(x,\tau)$ denote its inverse. Since $\chi\in C^1(\mathbb{R}\times[0,\tau^*))$ and $\chi_y>0$, the implicit function theorem gives $\xi\in C^1(\mathbb{R}\times[0,\tau^*))$, with
$$
\xi_x=\frac{1}{1+\tau u_0'(\xi)}
\qquad \mbox{and} \qquad
\xi_\tau=-\frac{u_0(\xi)}{1+\tau u_0'(\xi)} .
$$
Hence, in $\mathbb{R}\times(0,\tau^*)$, $N_\tau=\emptyset$ and $v(x,\tau)=(x-y_+(x,\tau))/\tau=(x-\xi(x,\tau))/\tau=u_0(\xi(x,\tau))$ by \eqref{X}. The function $v=u_0\circ\xi$ belongs to $C^1(\mathbb{R}\times[0,\tau^*))$, takes the value $u_0$ at $\tau=0$, and satisfies
$$
v_x=\frac{u_0'(\xi)}{1+\tau u_0'(\xi)}
\qquad \mbox{ and } \qquad
v_\tau=-\frac{u_0(\xi)u_0'(\xi)}{1+\tau u_0'(\xi)}=-vv_x .
$$
If $\tau^*\ge A_\infty$ we set $t^*:=\infty$. Otherwise $\kappa>0$, $\tau^*=1/\kappa<A_\infty$, and $t^*=\theta(1/\kappa)$ is the time defined in Theorem~\ref{thm:threshold}(ii). For $0\le t<t^*$ we have $A(t)<\tau^*$, and $u(x,t)=e^{-S(t)}v(x,A(t))$ is of class $C^1$, is given by \eqref{characteristics}, and is a classical (hence weak) solution of \eqref{main}, by the computation in Lemma~\ref{lem:viscous-reduction} with $\varepsilon=0$. Moreover
$$
u_x(x,t)=e^{-S(t)}\,\frac{u_0'(\xi)}{1+A(t)u_0'(\xi)},\qquad \xi=\xi(x,A(t)) .
$$
Since $1+A(t)u_0'(\xi)\ge1-\kappa A(t)>0$, the quotient $u_0'(\xi)/(1+A(t)u_0'(\xi))$ is an increasing function of $u_0'(\xi)\in[-\kappa,\infty)$. Moreover, $\inf u_0'=-\kappa$ when $\kappa>0$, and $\xi(\cdot,A(t))$ is onto $\mathbb{R}$. Therefore $u_x(x,t)\ge-\kappa e^{-S(t)}/(1-\kappa A(t))$, and, if $\kappa>0$, $\sup\limits_x(-u_x(x,t))=\kappa e^{-S(t)}/(1-\kappa A(t))$. For the upper bound, if $u_0'(\xi)\ge0$, the quotient is at most $u_0'(\xi)\le\|u_0'\|_{L^\infty}$, while if $u_0'(\xi)<0$, it is negative.

(i) If $\kappa A_\infty\le1$, then $\tau^*\ge A_\infty$, so $t^*=\infty$.

(ii) If $\kappa A_\infty>1$, then $t^*=\theta(1/\kappa)<\infty$. As $t\to(t^*)^-$, $A(t)$ increases to $1/\kappa$ and $e^{-S(t)}\to e^{-S(t^*)}>0$, which gives \eqref{blowup}.
\end{proof}

In both cases the classical solution coincides everywhere in its domain with the function \eqref{u}, since the minimizer is unique and equals $\xi(x,A(t))$. In particular, it is the entropy solution of Theorem~\ref{thm:uniqueness}, and for $t>t^*$ the solution is still given by \eqref{u}.

We now illustrate the role of $A_\infty$. In the first example, the classical threshold for constant damping is recovered.

\begin{example}[Constant coefficients]\label{ex:constant}
Let $\alpha\equiv\alpha_0>0$ and $\sigma\equiv\sigma_0>0$. Then $A(t)=\frac{\alpha_0}{\sigma_0}(1-e^{-\sigma_0t})$ and $A_\infty=\alpha_0/\sigma_0$. The solution remains classical for all times if and only if $\kappa\le\sigma_0/\alpha_0$. Otherwise the gradient blows up at $t^*=-\frac1{\sigma_0}\ln\big(1-\frac{\sigma_0}{\alpha_0\kappa}\big)$. This agrees with the elementary characteristic analysis of the damped Burgers equation. In the acoustic interpretation \eqref{gbe}, where $\alpha_0=1$, a constant $\sigma_0$ corresponds to $h(t)=h(0)e^{2\sigma_0t}$, that is, to an exponential horn, a duct whose cross-sectional area grows exponentially. In such a horn, waveforms whose compression $\kappa$ does not exceed $\sigma_0$ never break.
\end{example}

\begin{figure}[ht]
\centering
\begin{tikzpicture}[x=1.15cm,y=1.6cm,>=stealth]
\foreach \kap/\tmax/\sh in {0.8/2.6/0, 1.25/1.6094/7.4}{
\begin{scope}[xshift=\sh cm]
\draw[->] (-2.7,0)--(2.8,0) node[right] {$x$};
\draw[->] (-2.7,0.12)--(-2.7,2.85) node[above] {$t$};
\foreach \xi in {-2,-1.5,-1,-0.5,0,0.5,1,1.5,2}{
\draw[domain=0:\tmax,samples=50,variable=\t] plot ({(\xi)+(1-exp(-\t))*(-\kap*(\xi)/(1+(\xi)*(\xi)))},{\t});
}
\end{scope}
}
\node[below] at (0,-0.1) {$\kappa A_\infty\le1$};
\begin{scope}[xshift=7.4cm]
\draw[dashed] (-2.5,1.6094)--(2.5,1.6094);
\node[above] at (1.9,1.6094) {$t=t^*$};
\node[below] at (0,-0.1) {$\kappa A_\infty>1$};
\end{scope}
\end{tikzpicture}
\caption{Characteristics of \eqref{main} issued from smooth initial data. If $\kappa A_\infty\le1$ (left), they come to a halt before they can cross and the solution remains smooth for all times. If $\kappa A_\infty>1$ (right), they focus at the finite time $t^*$, where the gradient of the solution blows up.}
\label{fig:characteristics}
\end{figure}
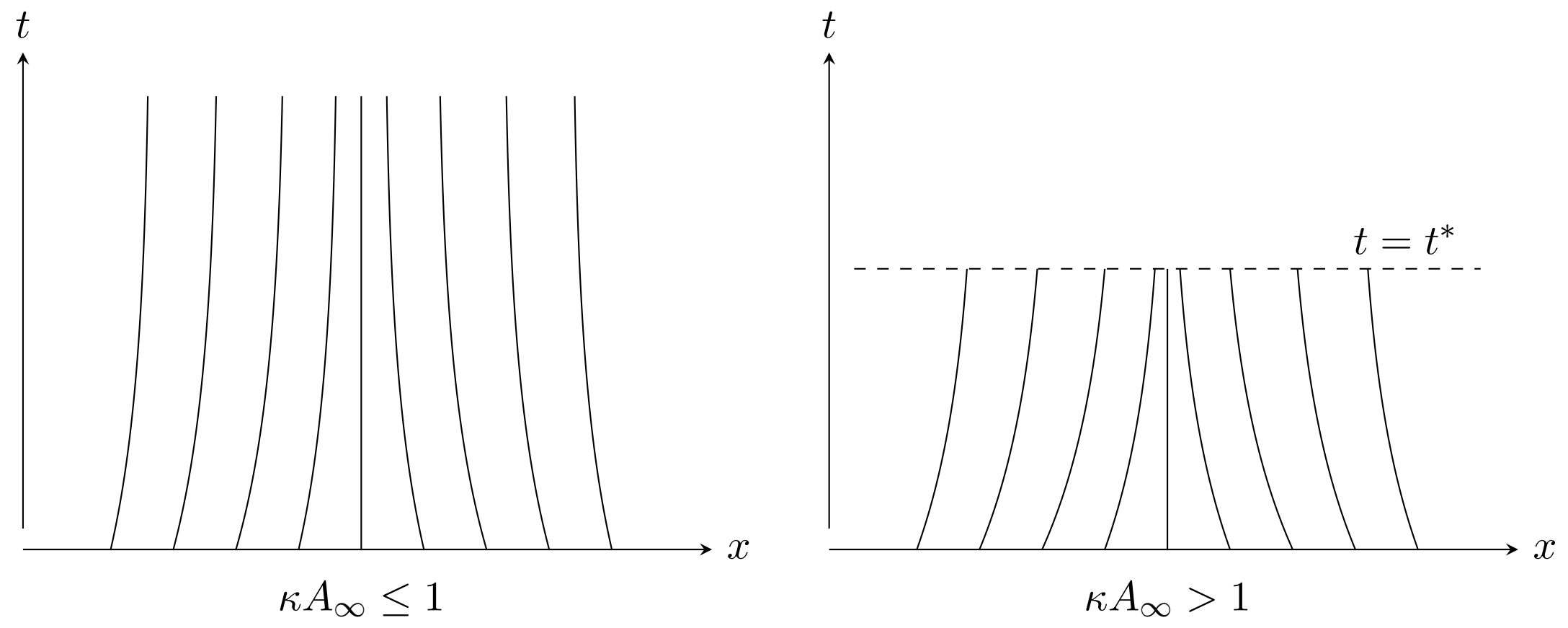

Figure~\ref{fig:characteristics} illustrates the two regimes of Theorem~\ref{thm:threshold} in the setting of Example~\ref{ex:constant} with $\alpha_0=\sigma_0=1$, so that $A(t)=1-e^{-t}$ and $A_\infty=1$. The initial datum is $u_0(\xi)=-\kappa\xi/(1+\xi^2)$, for which $\inf\limits_{\mathbb{R}}u_0'=u_0'(0)=-\kappa$, and the curves are the characteristics $x=\xi+A(t)u_0(\xi)$ of \eqref{characteristics} issued from $\xi=0,\pm0.5,\dots,\pm2$. The left panel corresponds to $\kappa=0.8$ and the right panel to $\kappa=1.25$, for which $t^*=\ln5$.

The next example corresponds to the generalized Burgers equation \eqref{gbe} for waves whose cross-sectional area grows like a power of the propagation distance, which includes cylindrical and spherical waves. For $t_0=1$, it is the inviscid nonplanar Burgers equation considered in \cite{RaoYadav2011}.

\begin{example}[Damping decaying like $1/t$ and nonplanar waves]\label{ex:mu}
Let $\alpha\equiv1$ and $\sigma(t)=m/(t_0+t)$, with $m>0$ and $t_0>0$. Then $e^{-S(t)}=\big(\frac{t_0}{t_0+t}\big)^{m}$ and
$$
A(t)=\begin{cases}\dfrac{t_0}{1-m}\Big[\Big(1+\dfrac{t}{t_0}\Big)^{1-m}-1\Big], & m\ne1,\\[10pt]
t_0\ln\Big(1+\dfrac{t}{t_0}\Big), & m=1,\end{cases}
\qquad
A_\infty=\begin{cases}\dfrac{t_0}{m-1}, & m>1,\\[8pt] \infty, & 0<m\le1.\end{cases}
$$
Hence for $m>1$ the solution remains classical for all times if and only if $\kappa\le(m-1)/t_0$, whereas for $0<m\le1$ the gradient blows up in finite time for every $u_0$ with $\kappa>0$, at
$$
t^*=t_0\Big[\Big(1+\frac{1-m}{\kappa t_0}\Big)^{1/(1-m)}-1\Big]\quad(m<1),
\qquad t^*=t_0\big(e^{1/(\kappa t_0)}-1\big)\quad(m=1).
$$
The value $m=1$ is thus critical. With $\sigma=h'/(2h)$ and $h(t)\propto(t_0+t)^j$, the cases $j=1$ and $j=2$ correspond to $m=1/2$ and $m=1$, respectively. Thus geometric spreading delays, but does not prevent, the breaking of cylindrical and spherical waves, and for spherical waves the breaking time is exponentially large with respect to $1/(\kappa t_0)$, in agreement with the classical theory of nonlinear spherical waves \cite{Whitham1974}.
\end{example}

In the absence of damping, the coefficient $\alpha$ alone may prevent shock formation.

\begin{example}[No damping, integrable coefficient $\alpha$]\label{ex:alpha}
Let $\sigma\equiv0$. Then $A_\infty=\int_0^\infty\alpha(t)\,dt$. For instance, if $\alpha(t)=(1+t)^{-2}$, then $A_\infty=1$ and the solution remains classical for all times if and only if $\kappa\le1$. Physically, the nonlinear transport switches off before characteristics have time to cross. In the acoustic interpretation, this is the second form of the generalized Burgers equation described in the Introduction, with $\alpha=(h(0)/h)^{1/2}$ and $h(t)=h(0)(1+t)^4$. Accordingly, the maximal displacement $A_\infty=1$ coincides with that of Example~\ref{ex:mu} with $m=2$ and $t_0=1$, because the two equations are related by the change of variable $u=(h(0)/h)^{1/2}v$.
\end{example}

Finally, when $\sigma$ and $\alpha$ are linked as in \cite{SuarezDelacruz2026}, their effects cancel exactly.

\begin{example}[Compensated damping]\label{ex:compensated}
If $\alpha\in C^1([0,\infty))$ and $\sigma=\alpha'/\alpha$, then $a\equiv\alpha(0)$ and $A(t)=\alpha(0)t$, so $A_\infty=\infty$. Hence the gradient blows up at $t^*=1/(\alpha(0)\kappa)$ for every $u_0$ with $\kappa>0$, regardless of the size of the damping.
\end{example}

\section{Solution of the Riemann problem}\label{sec:riemann}

For Riemann data the variational problem can be solved by hand, which provides a transparent illustration of the general theory and allows a comparison with the solutions obtained in \cite{DelacruzLuWang2024} by a different nonautonomous regularization. The computation also shows that, on the shock curve, the minimizer is not unique and the two minimizers lie on opposite sides of the origin.

\begin{proof}[Proof of Theorem~\ref{thm:riemann}]
Since $u_0\in L^\infty(\mathbb{R})$, by Theorem~\ref{thm:uniqueness} it suffices to compute $v=(x-y_+)/\tau$ and to use $u(x,t)=e^{-S(t)}v(x,A(t))$. Here $U_0(y)=u_ly$ for $y<0$ and $U_0(y)=u_ry$ for $y\ge0$. On $(-\infty,0]$ the function $G(x,\cdot,\tau)$ is a quadratic polynomial with critical point $y_l:=x-u_l\tau$, and on $[0,\infty)$ it is a quadratic polynomial with critical point $y_r:=x-u_r\tau$. A direct computation gives
\begin{equation}\label{Gl-Gr}
G(x,y_l,\tau)=u_lx-\frac{u_l^2}{2}\tau,\qquad G(x,y_r,\tau)=u_rx-\frac{u_r^2}{2}\tau,\qquad
G(x,0,\tau)=\frac{x^2}{2\tau}.
\end{equation}

\emph{Case $u_l<u_r$.} Then $U_0$ is convex, so $G(x,\cdot,\tau)$ is strictly convex and has a unique minimizer. If $x<u_l\tau$, then $y_l<0$ is a critical point of $G(x,\cdot,\tau)$, hence the minimizer, and $v=u_l$. If $x>u_r\tau$, then $y_r>0$ is the minimizer and $v=u_r$. If $u_l\tau\le x\le u_r\tau$, then the left and right derivatives of $G(x,\cdot,\tau)$ at $y=0$ are $-x/\tau+u_l\le0$ and $-x/\tau+u_r\ge0$, respectively, so the minimizer is $0$ and $v=x/\tau$. This gives \eqref{rarefaction}.

\emph{Case $u_l>u_r$.} The minimum of $G(x,\cdot,\tau)$ over $(-\infty,0]$ is attained at $\min\{y_l,0\}$ and over $[0,\infty)$ at $\max\{y_r,0\}$. If $x\le u_r\tau$, then $y_l<y_r\le0$, and $G(x,y_l,\tau)<x^2/(2\tau)=G(x,0,\tau)$ because $y_l\ne0$. Since in this case the minimum over $[0,\infty)$ is $G(x,0,\tau)$, the unique minimizer is $y_l$. Symmetrically, if $x\ge u_l\tau$, the unique minimizer is $y_r$. If $u_r\tau<x<u_l\tau$, both $y_l<0<y_r$ are admissible, both values in \eqref{Gl-Gr} are smaller than $G(x,0,\tau)$, and
$$
G(x,y_l,\tau)-G(x,y_r,\tau)=(u_l-u_r)\Big(x-\frac{u_l+u_r}{2}\tau\Big).
$$
Hence the unique minimizer is $y_l$ for $x<\frac{u_l+u_r}{2}\tau$ and $y_r$ for $x>\frac{u_l+u_r}{2}\tau$, while for $x=\frac{u_l+u_r}{2}\tau$ there are exactly two minimizers, $y_l=\frac{u_r-u_l}{2}\tau<0$ and $y_r=\frac{u_l-u_r}{2}\tau>0$. Therefore $v=u_l$ for $x<\frac{u_l+u_r}{2}\tau$ and $v=u_r$ for $x>\frac{u_l+u_r}{2}\tau$, which gives \eqref{shock}.
\end{proof}

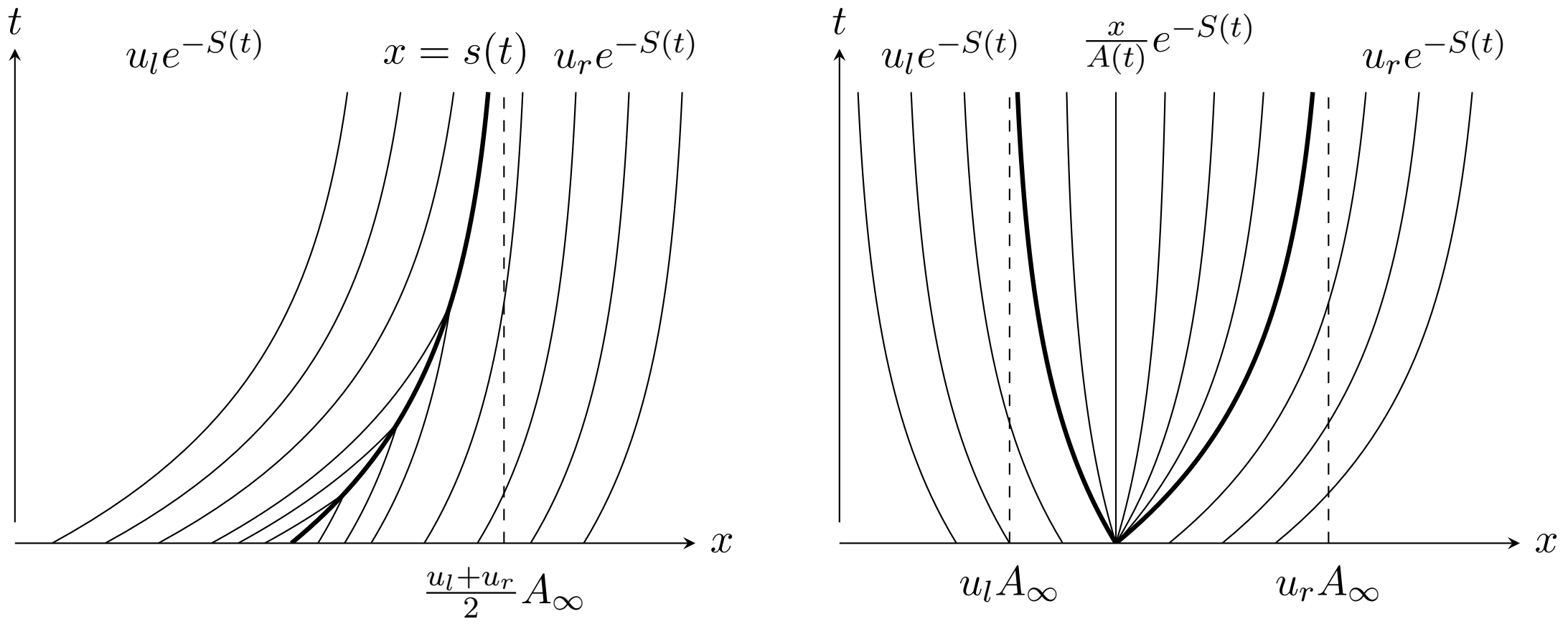
\begin{figure}[!ht]
\centering
\begin{tikzpicture}[x=0.98cm,y=1.6cm,>=stealth]
\draw[->] (-2.6,0)--(3.8,0) node[right] {$x$};
\draw[->] (-2.6,0.12)--(-2.6,2.85) node[above] {$t$};
\begin{scope}
\clip (-2.5,0) rectangle (3.7,2.6);
\foreach \xi in {-2.25,-1.75,-1.25}{\draw[thin,domain=0:2.6,samples=40,variable=\t] plot ({\xi+3*(1-exp(-\t))},{\t});}
\foreach \xi/\te in {-0.75/1.3863,-0.5/0.6931,-0.25/0.2877}{\draw[thin,domain=0:\te,samples=30,variable=\t] plot ({\xi+3*(1-exp(-\t))},{\t});}
\foreach \xi/\te in {0.25/0.2877,0.5/0.6931,0.75/1.3863}{\draw[thin,domain=0:\te,samples=30,variable=\t] plot ({\xi+(1-exp(-\t))},{\t});}
\foreach \xi in {1.25,1.75,2.25,2.75}{\draw[thin,domain=0:2.6,samples=40,variable=\t] plot ({\xi+(1-exp(-\t))},{\t});}
\draw[very thick,domain=0:2.6,samples=60,variable=\t] plot ({2*(1-exp(-\t))},{\t});
\draw[dashed] (2,0)--(2,2.6);
\end{scope}
\node[below] at (2,-0.05) {$\frac{u_l+u_r}{2}A_\infty$};
\node[above] at (-0.9,2.62) {$u_le^{-S(t)}$};
\node[above] at (1.55,2.62) {$x=s(t)$};
\node[above] at (3.15,2.62) {$u_re^{-S(t)}$};
\begin{scope}[xshift=7.6cm]
\draw[->] (-2.6,0)--(3.8,0) node[right] {$x$};
\draw[->] (-2.6,0.12)--(-2.6,2.85) node[above] {$t$};
\begin{scope}
\clip (-2.5,0) rectangle (3.7,2.6);
\foreach \xi in {-1.5,-1,-0.5}{\draw[thin,domain=0:2.6,samples=40,variable=\t] plot ({\xi-(1-exp(-\t))},{\t});}
\foreach \xi in {0.5,1,1.5}{\draw[thin,domain=0:2.6,samples=40,variable=\t] plot ({\xi+2*(1-exp(-\t))},{\t});}
\foreach \c in {-0.5,0,0.5,1,1.5}{\draw[thin,domain=0:2.6,samples=40,variable=\t] plot ({\c*(1-exp(-\t))},{\t});}
\foreach \c in {-1,2}{\draw[very thick,domain=0:2.6,samples=60,variable=\t] plot ({\c*(1-exp(-\t))},{\t});}
\draw[dashed] (-1,0)--(-1,2.6);
\draw[dashed] (2,0)--(2,2.6);
\end{scope}
\node[below] at (-1,-0.05) {$u_lA_\infty$};
\node[below] at (2,-0.05) {$u_rA_\infty$};
\node[above] at (-1.55,2.62) {$u_le^{-S(t)}$};
\node[above] at (0.5,2.62) {$\frac{x}{A(t)}e^{-S(t)}$};
\node[above] at (3.0,2.62) {$u_re^{-S(t)}$};
\end{scope}
\end{tikzpicture}
\caption{Solution of the Riemann problem when $A_\infty<\infty$. The shock curve (left) and the edges of the rarefaction fan (right) approach vertical asymptotes, so that both waves remain confined to a bounded region of space while their amplitudes decay.}
\label{fig:riemann}
\end{figure}

In the shock case, $s'(t)=\frac{u_l+u_r}{2}a(t)=\alpha(t)\frac{u(s(t)-,t)+u(s(t)+,t)}{2}$, which is the Rankine--Hugoniot condition for \eqref{main}, and $u(s(t)-,t)>u(s(t)+,t)$. In the rarefaction case the fan occupies the region $u_lA(t)\le x\le u_rA(t)$, and the solution is self-similar in $x/t$ only when $A$ is linear. If $A_\infty<\infty$, the waves are confined, as shown in Figure~\ref{fig:riemann} for $\alpha\equiv1$ and $\sigma\equiv1$, where $A(t)=1-e^{-t}$ and $A_\infty=1$, with $(u_l,u_r)=(3,1)$ in the shock case and $(u_l,u_r)=(-1,2)$ in the rarefaction case. The shock position converges to $\frac{u_l+u_r}{2}A_\infty$, and the rarefaction fan converges to the bounded interval $[u_lA_\infty,u_rA_\infty]$, while the amplitudes decay like $e^{-S(t)}$. The solutions \eqref{rarefaction} and \eqref{shock} coincide with those obtained in \cite[Theorem~10]{DelacruzLuWang2024} for $k=\lambda=1$ by means of a different nonautonomous viscous regularization, and Theorem~\ref{thm:riemann} additionally yields their uniqueness in the class of entropy solutions.

\section*{Acknowledgments}
The authors thank the Universidad Pedag\'ogica y Tecnol\'ogica de Colombia for the academic environment and the research time that made this work possible.

\section*{Author Contributions}
Both authors contributed equally to the conception of this work, the mathematical analysis, and the writing of the manuscript, and both approved its final version.

\section*{Financial Disclosure}
None reported.

\section*{Conflicts of Interest}
The authors declare that there are no potential conflicts of interest.

\end{document}